\documentclass[12pt, leqno]{article}
\usepackage{amsmath,amsfonts}
\usepackage{amssymb}
\usepackage{amsthm}
\usepackage{times,txfonts}

\usepackage{accents}

\usepackage{pict2e}

\newtheorem{thm}{Theorem}[section]
\newtheorem{prop}[thm]{Proposition}

\newtheorem{lemma}[thm]{Lemma}
\newtheorem{dfn}[thm]{Definition}

\newtheorem{remark}[thm]{Remark}
\newtheorem{example}[thm]{Example}

\numberwithin{equation}{section}

\def\pf{\noindent{\it Proof.} \ }

\def\qed{\hfill $\square$}

\def\id{{\rm id}}

\title{UC hierarchy and monodromy preserving deformation}

\author{Teruhisa Tsuda}
\date{January 21, 2010 (Revised: January 31, 2012)}

\begin{document}
\maketitle

\begin{abstract}
The UC hierarchy is an extension of the KP hierarchy,
which possesses not only an infinite set of positive time evolutions 
but also that of negative ones.
Through a similarity reduction
we derive from the UC hierarchy
a class of the Schlesinger systems
including the Garnier system and
the sixth Painlev\'e equation,
which describes the monodromy preserving deformations of Fuchsian linear differential equations with certain spectral types.
We also present a unified formulation of 
the above Schlesinger systems
as a canonical Hamiltonian system whose Hamiltonian functions are polynomials in the canonical variables.
\end{abstract}

\renewcommand{\thefootnote}{\fnsymbol{footnote}}

\section{Introduction}

This work is aimed to present
a certain connection between 
infinite-dimensional integrable systems of soliton type and 
finite-dimensional integrable systems of isomonodromic type.
The KP (Kadomtsev--Petviashvili) hierarchy is, undoubtedly, the most basic one
among the former
and is
a series of nonlinear partial differential equations in infinitely many 
independent variables
${\boldsymbol x}=(x_1,x_2,x_3,\ldots)$
that are consistent with each other.
It literally includes 
as the first nontrivial member
the KP  equation
\begin{equation}  \label{eq:kp_eq} 
\frac{3}{4} \frac{\partial^2 f }{\partial {x_2}^2}
=  \frac{\partial }{\partial x_1}
\left(
\frac{\partial f }{\partial x_3} -\frac{3}{2} f \frac{\partial f }{\partial x_1}
-\frac{1}{4}  \frac{\partial^3 f }{\partial {x_1}^3}
\right),
\end{equation}
which is a typical soliton equation.
If we count the degree of variables as $\deg x_n=n$ 
and $\deg f =-2$,
then both sides of (\ref{eq:kp_eq}) are 
equally homogeneous (of degree $-6$)
as differential polynomials in $f$
with respect to $x_n$.
Every equation of the KP hierarchy is known to be homogeneous, in fact.
In this sense
we may say that the KP hierarchy forms a 
homogeneous integrable system
equipped with an infinite set of
time evolutions of positive degree.
The {\it UC} ({\it Universal Character}) {\it hierarchy},
introduced in \cite{tsu04},
 is an infinite-dimensional integrable system
which naturally generalizes the KP hierarchy by taking into account 
the negative time evolutions besides the positive ones 
while keeping its homogeneity.
The independent variables of the UC hierarchy
consist of two sets of infinitely many variables
${\boldsymbol x}$ and ${\boldsymbol y}=(y_1,y_2,y_3, \ldots)$
with their degrees given as $\deg x_n=n$ and $\deg y_n=-n$.
In this paper
we show  that a similarity reduction of the UC hierarchy yields 
a broad class of the Schlesinger systems
including the Garnier system and the sixth Painlev\'e equation,
which describes the monodromy preserving deformations of 
Fuchsian linear differential equations
with certain spectral types.

We begin by recalling the definition of the UC hierarchy.
Let us introduce the commuting pair of 
linear differential operators
(called the {\it vertex operators})
\begin{subequations}  \label{subeq:vo}
\begin{align}
X^\pm(z)&= \sum_{n \in {\mathbb Z}} X_n^\pm z^n =e^{\pm \xi({\boldsymbol x}-\widetilde{\partial}_{\boldsymbol y},z)} e^{\mp \xi(\widetilde{\partial}_{\boldsymbol x},z^{-1})}, \\
Y^\pm(w)&= \sum_{n \in {\mathbb Z}} Y_n^\pm w^n =e^{\pm \xi({\boldsymbol y}-\widetilde{\partial}_{\boldsymbol x},w)} e^{\mp \xi(\widetilde{\partial}_{\boldsymbol y},w^{-1})}, 
\end{align}
\end{subequations}
where we have used the notations
\[\xi({\boldsymbol x},z)=\sum_{n=1}^\infty x_n z^n
\quad
\text{and}
\quad
 \widetilde{\partial}_{\boldsymbol x}=\left( \frac{\partial}{\partial x_1},\frac{1}{2}\frac{\partial}{\partial x_2},
\frac{1}{3}\frac{\partial}{\partial x_3},\ldots \right).
\]

\begin{dfn}  \rm   
For an unknown function $\tau=\tau({\boldsymbol x},{\boldsymbol y})$,
the simultaneous bilinear equation
\begin{equation} \label{eq:UCH}
\sum_{i+j=-1}X_i^{-} \tau \otimes X_j^{+} \tau =
\sum_{i+j=-1}Y_i^{-} \tau \otimes Y_j^{+} \tau =0
\end{equation}
is called  the {\it UC hierarchy}.
\end{dfn}

The UC hierarchy is homogeneous indeed
as it has the following scaling symmetry:
if $\tau$ is a solution 
 of (\ref{eq:UCH})
 then so 
is
$\tau(c x_1,c^2 x_2,\ldots,c^{-1} y_1,c^{-2} y_2,\ldots)$
for any $c \in {\mathbb C}^\times$.
The UC hierarchy is regarded as an extension of the KP hierarchy.
 If $\tau$ does not depend on ${\boldsymbol y}$, then the latter equality 
of (\ref{eq:UCH})
trivially holds
 and the former reduces to the bilinear expression of the KP hierarchy,
 which is due to Date--Jimbo--Kashiwara--Miwa
 (see \cite{kac90, mjd00});
 for reference the variable transformation toward 
 the original KP equation, (\ref{eq:kp_eq}),
 is given by 
 $f=2 (\partial/\partial x_1)^2 \log \tau$.
 We always require the solution
$\tau=\tau({\boldsymbol x},{\boldsymbol y})$,
called the {\it $\tau$-function},
 to be an entire function 
with respect to each independent variable.
Note that $\tau$-functions are distinguished up to multiplication by constants,
as can be seen from (\ref{eq:UCH}). 
Concerning the UC hierarchy there is a counterpart of the Sato theory about the KP hierarchy; cf. \cite{sat}.
That is,
the totality of solutions 
of the UC hierarchy
forms a direct product of two Sato Grassmannians and the action of its transformation group can be realized by means of the vertex operators.
For details to \cite{tsu04}.
Of particular interest is its homogeneous polynomial solution,
which is a fixed solution with respect to the above scaling symmetry.

Let $\lambda=(\lambda_1,\lambda_2, \ldots,\lambda_{\ell})$ and 
$\mu=(\mu_1,\mu_2, \ldots, \mu_{\ell'})$
be a pair of partitions.
Consider the following determinant of {\it twisted} Jacobi--Trudi type:
\begin{equation}  \label{eq:twist}
S_{[\lambda,\mu ]}
= \det 
\left(
  \begin{array}{ll}
 \widetilde{h}_{\mu_{\ell'-i+1}  +i - j },  &  i \leq \ell'  \\
 h_{\lambda_{i-\ell'}-i+j},     &   i > \ell'+ 1     \\
  \end{array}
\right)_{1 \leq i,j \leq \ell+\ell'},
\end{equation}
where
$h_n=h_n({\boldsymbol x})$ $(n \in {\mathbb Z})$
is a polynomial in only ${\boldsymbol x}$ 
and is
defined by the generating function
\[e^{\xi({\boldsymbol x},z)}=\sum_{n \in {\mathbb Z}} h_n z^n,
\]
and $\widetilde{h}_n=\widetilde{h}_n({\boldsymbol y})$ 
is exactly the same as $h_n$
except replacing ${\boldsymbol x}$ with ${\boldsymbol y}$.
If $\mu=\emptyset$ then (\ref{eq:twist}) reduces to  the (usual) Jacobi--Trudi formula: $S_\lambda=S_{[\lambda,\emptyset]}=\det  (h_{\lambda_{i}-i+j})$, which defines the Schur function 
$S_\lambda=S_\lambda({\boldsymbol x})$.
The polynomial $S_{[\lambda,\mu ]}=S_{[\lambda,\mu ]}({\boldsymbol x},{\boldsymbol y})$
is called the {\it universal character} and was
originally
 introduced by Koike \cite{koi89}
 in the study of classical groups.
It is easy to see 
that 
$S_{[\lambda,\mu ]}$
becomes a homogeneous polynomial whose degree equals 
the difference $|\lambda|-|\mu|$,
where the sum 
 $|{\lambda}|=\lambda_1+\cdots+\lambda_\ell$ 
denotes
the weight of a partition $\lambda$.
A few examples are
$S_{[\emptyset,\emptyset]}=1$, $S_{[(1),\emptyset]}=x_1$, 
$S_{[(1),(1)]}=x_1y_1-1$, $S_{[(2,1),(1)]}=y_1\left({x_1}^3/3-x_3\right)-{x_1}^2$, 
etc.
Remarkably,
the set of homogeneous polynomial solutions of the UC hierarchy,
(\ref{eq:UCH}),
coincides with
that of the universal characters
$\{S_{[\lambda,\mu]}({\boldsymbol x},{\boldsymbol y})\}_{\lambda,\mu: \text{\rm partitions}}$.

By considering general
homogeneous solutions of the UC hierarchy
that are not necessarily polynomials,
we can find a link to the theory of monodromy preserving deformations.
Let us explain it in more detail.
First we derive from the original one, 
(\ref{eq:UCH}),
similar bilinear equations among solutions
generated by successive application of vertex operators.
Let
$\tau_{m,n}=\tau_{m,n}({\boldsymbol x},{\boldsymbol y})$
denote such a sequence of
solutions of the UC hierarchy.
A typical example of the bilinear equations 
is 
\[
\tau_{m,n} \otimes \tau_{m+1,n+1}
= \sum_{i+j=0} X_i^{-} \tau_{m+1,n} \otimes X_j^{+} \tau_{m,n+1}.
\]
Next we impose on the sequence $\tau_{m,n}$ of solutions
homogeneity
\begin{equation} \label{eq:homog}
E \tau_{m,n}=d_{m,n} \tau_{m,n} \quad
(d_{m,n} \in {\mathbb C})
\end{equation}
and periodicity
\begin{equation} \label{eq:period}
\tau_{m+L,n}=\tau_{m,n+L}=\tau_{m,n}
\end{equation}
for an integer $L(\geq 2)$ fixed. 
Here we have used the {\it Euler operator}
\[
E=
 \sum_{n=1}^\infty 
\left(n x_n \frac{\partial}{\partial x_n} - n y_n \frac{\partial}{\partial y_n}
\right),
\]
which is a linear differential operator measuring the degree of a homogeneous function; for instance, $E S_{[\lambda,\mu]}=(|\lambda|-|\mu|)  S_{[\lambda,\mu]}$.
Finally we substitute into each 
$x_n$ and $y_n$ the `power sum'
of new independent variables ${\boldsymbol t}=(t_0,t_1,\ldots,t_N)$
as
\begin{equation} \label{eq:power}
x_n=\frac{1}{n}\sum_{i=0}^N \theta_i {t_i}^n \quad \text{and} \quad
y_n=\frac{1}{n}\sum_{i=0}^N\theta_i {t_i}^{-n}
\quad 
(n=1,2,\ldots)
\end{equation}
where $\theta_i \in {\mathbb C}$ are constant parameters.
In view of the homogeneity (\ref{eq:homog}),
we may take $t_0=1$ without loss of generality. 
Under the reduction conditions 
(\ref{eq:homog}), (\ref{eq:period}), and (\ref{eq:power}),
the UC hierarchy yields
a system of nonlinear partial differential equations
in $N$ variables,
hereafter denoted by ${\cal G}_{L,N}$, 
whose phase space 
is essentially of $2N(L-1)$ dimension
(see Theorem~\ref{thm:gar}).
To sum up the above procedure,
we say that 
${\cal G}_{L,N}$ 
is a {\it similarity reduction}
of the UC hierarchy.
The system ${\cal G}_{L,N}$ 
is a finite-dimensional integrable system of isomonodromic type.
For instance
${\cal G}_{2,N}$ 
corresponds to the 
Garnier system in $N$ variables
and 
${\cal G}_{2,1}$,
the first nontrivial case, 
does the sixth Painlev\'e equation.
From the viewpoint of the UC hierarchy
we can clearly understand various aspects of ${\cal G}_{L,N}$,
e.g., 
Hirota bilinear relations for $\tau$-functions,
Weyl group symmetries, 
and algebraic solutions expressed in terms of the universal character.

As analogous to the case of the KP hierarchy,
the UC hierarchy (\ref{eq:UCH})
generates
the linear equations for unknown functions
(called the {\it wave functions})
\[
\psi_{m,n}=
\psi_{m,n}({\boldsymbol x},{\boldsymbol y},k)
= \frac{ \tau_{m,n-1}({\boldsymbol x}-[k^{-1}],{\boldsymbol y}-[k])  }{ \tau_{m,n}({\boldsymbol x},{\boldsymbol y}) }e^{\xi({\boldsymbol x},k)},
\]
where $[k]=(k,k^2/2,k^3/3,\ldots)$.
Through the reduction procedure
they induce an auxiliary system of linear differential equations;
one of which is
a Fuchsian system
of rank $L$
 in
the {\it spectral variable}
$z=k^L$
with $N+3$ poles on the Riemann sphere,
and the others govern its monodromy preserving deformations.
The nonlinear system
${\cal G}_{L,N}$
can be reformulated as 
a compatibility condition of this auxiliary linear system (Lax formalism).
Remark here that the compatibility itself is {\it a priori} established
because all the linear equations originate from the same 
bilinear equation (\ref{eq:UCH}).

The {\it spectral type} of the 
Fuchsian system
under consideration
is given by 
 the $(N+3)$-tuple 
\[
 \underbrace{(L-1,1), \ldots, (L-1,1)}_{N+1},
(1,1,\ldots,1),(1,1,\ldots,1)
\]
of partitions of $L$,
which indicates how the characteristic exponents 
overlap
at each of the $N+3$ singularities.
Thus
we conclude that
${\cal G}_{L,N}$ is equivalent to
 a particular case of the Schlesinger systems
 specified by this spectral type.
We also present 
a unified description of ${\cal G}_{L,N}$ 
for any $L$ and $N$
as a canonical Hamiltonian system,
denoted by ${\cal H}_{L,N}$,
whose Hamiltonian functions are
polynomials in the canonical variables
(see Theorem~\ref{thm:caneq}).

\begin{remark}\rm
In \cite{tsu10}
we found particular solutions of ${\cal G}_{L,N}$ (or ${\cal H}_{L,N}$)
expressed in terms of a certain generalization of Gau\ss's hypergeometric function.
\end{remark}

In the next section
we derive some difference (and differential) equations
from the UC hierarchy
as a preliminary.
In Sect.~\ref{sect:hom}, 
we construct a sequence of homogeneous solutions
of the UC hierarchy and present
its Weyl group symmetry of type $A$.
In Sect.~\ref{sect:sim},
we consider a similarity reduction of the UC hierarchy 
by requiring its solutions to satisfy the homogeneity and periodicity.
As a result we obtain a nonlinear system
${\cal G}_{L,N}$
of partial differential equations,
which provides an extension of   
both the Garnier system and the sixth Painlev\'e equation.
The universal characters $S_{[\lambda,\mu]}$
are homogeneous solutions of 
the UC hierarchy and 
thereby consistent with the similarity reduction.
Hence, as described in Sect.~\ref{sect:ratsol},
it is immediate to obtain particular
solutions of ${\cal G}_{L,N}$ expressed in terms of $S_{[\lambda,\mu]}$.
The subject of Sect.~\ref{sect:lax} is the
Lax formalism of 
the systems 
${\cal G}_{L,N}$,
which reveals that they 
constitute a class of the Schlesinger systems.
We show that the auxiliary linear problem of ${\cal G}_{L,N}$
arises naturally from the linear equations satisfied by 
the wave functions of the UC hierarchy.
In Sect.~\ref{sect:ham},
we transform ${\cal G}_{L,N}$ into
the canonical Hamiltonian  
system ${\cal H}_{L,N}$
with polynomial Hamiltonian functions.
Section~\ref{sect:weyl} 
is devoted to 
the birational symmetries.
We observe that the Weyl group actions, 
discussed in Sect.~\ref{sect:hom}, 
give rise to birational canonical transformations
of ${\cal H}_{L,N}$.
In the appendix we briefly indicate a relationship between 
our polynomial Hamiltonian structure 
and that
given by Kimura and Okamoto \cite{ko84}
for the Garnier system, i.e.,
the case where $L=2$.

\section{Method for generating a `closed' functional equation}
\label{sect:fe}

Unlike in the case of the KP hierarchy, 
every differential equation 
of the UC hierarchy
with respect to
the original variables
${\boldsymbol x}$ and ${\boldsymbol y}$
 is of infinite order. 
In this section
we show how to overcome this difficulty, i.e., 
a method for generating 
a `closed' functional equation from the UC hierarchy;
cf. \cite{djm82}.

We first recall that if $\tau=\tau({\boldsymbol x},{\boldsymbol y})$ 
is a solution of (\ref{eq:UCH}) then so are $X^+(a)\tau$ and $Y^+(b)\tau$
for any
$a, b \in {\mathbb C}^\times$.
With this fact in mind,
let us take our interest in 
bilinear equations 
for a sequence 
of solutions generated 
by successive application of vertex operators.
Suppose 
$\tau_{0,0}=\tau({\boldsymbol x},{\boldsymbol y})$
to be a solution of the UC hierarchy, (\ref{eq:UCH}).
Define a sequence $ \tau_{m,n}$ of solutions by
\begin{equation} \label{eq:seq}
\tau_{m,n}= \prod_{i=0}^{m-1} X^+(a_i) 
 \prod_{j=0}^{n-1} Y^+(b_j) \tau_{0,0},
 \end{equation}
where we write as
\[
\prod_{i=0}^{m-1} X^+(a_i) =X^+(a_{m-1})\cdots X^+(a_{1})X^+(a_{0}).
\]
Then we can derive
similar bilinear equations from the UC hierarchy,
the original one (\ref{eq:UCH}).

\begin{lemma}  \label{lemma:mUC}
For integers $m,n \geq 0$,
it holds that 
\begin{align} 
&
\sum_{i+j=-m-1}X_i^{-} \tau_{0,0} \otimes X_j^{+} \tau_{m,n} =
\sum_{i+j=-n-1}Y_i^{-} \tau_{0,0} \otimes Y_j^{+} \tau_{m,n} =0,
\label{eq:mUC1}
\\
&
\tau_{0,0} \otimes \tau_{1,n}
- \sum_{i+j=0} X_i^{-} \tau_{1,0} \otimes X_j^{+} \tau_{0,n}
=\sum_{i+j=-n-1}Y_i^{-} \tau_{1,0} \otimes Y_j^{+} \tau_{0,n} =0,
\label{eq:mUC2}
\\
&
\sum_{i+j=-m-1}X_i^{-} \tau_{0,1} \otimes X_j^{+} \tau_{m,0} 
=
\tau_{0,0} \otimes \tau_{m,1}
- \sum_{i+j=0} Y_i^{-} \tau_{0,1} \otimes Y_j^{+} \tau_{m,0}=0.
\label{eq:mUC3}
\end{align}
\end{lemma}

\pf
Notice that 
 the operators $X_i^\pm$ ($i \in {\mathbb Z}$)
satisfy the {\it fermionic} relations:
$X_i^\pm X_j^\pm+X_{j-1}^\pm X_{i+1}^\pm =0$ and 
$X_i^+X_j^- +X_{j+1}^-X_{i-1}^+ =\delta_{i+j,0}$.
The same relations hold also for  $Y_i^\pm$.
Moreover, $X_i^\pm$  and $Y_j^\pm$ mutually commute. 
See \cite{tsu04}.
By virtue of the above relations,
applying 
$1 \otimes\prod_{i=0}^{m-1} X^+(a_i)\prod_{j=0}^{n-1} Y^+(b_j)$,
$X^+(a_0) \otimes \prod_{j=0}^{n-1} Y^+(b_j)$,
and
$Y^+(b_0) \otimes\prod_{i=0}^{m-1} X^+(a_i)$
to (\ref{eq:UCH}),
we obtain
(\ref{eq:mUC1}), (\ref{eq:mUC2}),
and (\ref{eq:mUC3}),
respectively.
\qed 
\\

We shall look closely at (\ref{eq:mUC1}),
which corresponds to the original UC hierarchy
(\ref{eq:UCH}) when $m=n=0$.
It can be rewritten
equivalently into
\begin{subequations}  \label{subeq:mn_uc}
\begin{align}
 &\frac{1}{2 \pi \sqrt{-1}} \oint 
 z^m e^{\xi({\boldsymbol x}-{\boldsymbol x'},z)} {\rm d}z \,
\tau_{0,0} ({\boldsymbol x'}+[z^{-1}],{\boldsymbol y'}+[z]) 
\tau_{m,n}({\boldsymbol x}-[z^{-1}],{\boldsymbol y}-[z])
=0,  \label{eq:mn_uc_1}
\\
&\frac{1}{2 \pi \sqrt{-1}} \oint 
w^n e^{\xi({\boldsymbol y}-{\boldsymbol y'},w)} {\rm d}w \,
\tau_{0,0}({\boldsymbol x'}+[w],{\boldsymbol y'}+[w^{-1}]) 
\tau_{m,n}({\boldsymbol x}-[w],{\boldsymbol y}-[w^{-1}])
=0   \label{eq:mn_uc_2}
\end{align}
\end{subequations}
with 
${\boldsymbol x}$, ${\boldsymbol y}$, ${\boldsymbol x'}$, and ${\boldsymbol y'}$
being arbitrary parameters,
where
$\oint \frac{{\rm d} z}{2 \pi \sqrt{-1}}  $
means taking the coefficient of $1/z$ of
the integrand as a (formal) 
Laurent series expansion in $z$.
If we try to write down a differential equation {\it naively} 
after the case of the KP hierarchy,
namely if we consider the Taylor series expansion of 
 (\ref{eq:mn_uc_1})
 around 
 $\{{\boldsymbol x'}={\boldsymbol x},{\boldsymbol y'}={\boldsymbol y}\}$,
then we have 
an infinite set of differential equations of infinite order;
see \cite{tsu04}. 
This result reflects the fact
that the integrand of 
 (\ref{eq:mn_uc_1})
under the substitution 
${\boldsymbol x'}={\boldsymbol x}$
and
 ${\boldsymbol y'}={\boldsymbol y}$
may admit
an essential singularity not only at $z=0$
but also at $z=\infty$.
However,
we can construct a functional equation
in a closed expression
by taking an appropriate choice of parameters
${\boldsymbol x}$, ${\boldsymbol y}$, ${\boldsymbol x'}$, and ${\boldsymbol y'}$
instead.

\begin{figure}
\begin{center}
\begin{picture}(120,100)
\thicklines

\put(0,90){$C_2$}
\put(40,50){\circle{80}}
\put(0,55){\vector(0,1){0}}
\put(80,45){\vector(0,-1){0}}
\put(15,57){$0$}
\put(16,48){\line(1,1){4}}
\put(16,52){\line(1,-1){4}}

\put(102,90){$C_1$}
\put(75,50){\circle{80}}
\put(35,55){\vector(0,1){0}}
\put(115,45){\vector(0,-1){0}}

\put(92,57){$\infty$}
\put(95,48){\line(1,1){4}}
\put(95,52){\line(1,-1){4}}

\put(58,60){${t_i}^{-1}$}

\put(48,58){\line(1,1){4}}
\put(48,62){\line(1,-1){4}}
\put(48,48){\line(1,1){4}}
\put(48,52){\line(1,-1){4}}
\put(49,32){$\vdots$}

\put(46,-5){($i \in I$)}
\end{picture}
\end{center}
\caption{Contours of integration and singularities of $F(z)$.}
\label{fig:cont}
\end{figure}
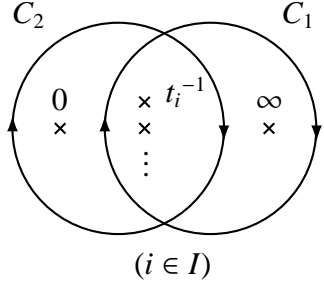

Let $I, J \subset {\mathbb Z}$ be a disjoint pair of finite indexing sets.
By specializing the parameters in (\ref{subeq:mn_uc}) 
as
\[
{\boldsymbol x'} ={\boldsymbol x}-\sum_{j \in I}[t_j]+ \sum_{j \in J} [t_j],
\quad
 {\boldsymbol y'} ={\boldsymbol y}-\sum_{j \in I}[{t_j}^{-1}]+ \sum_{j \in J} [{t_j}^{-1}],
\]
we obtain
\begin{align*}
\Omega_1
&:= 
 z^m e^{\xi({\boldsymbol x}-{\boldsymbol x'},z)} {\rm d}z
 =z^m\frac{\prod_{j \in J} (1-t_j z)}{ \prod_{j \in I} (1-t_j z)} {\rm d}z
 \quad (|t_j z| <1),
\\
\Omega_2
&:= 
 w^n e^{\xi({\boldsymbol y}-{\boldsymbol y'},w)} {\rm d}w
 =w^n\frac{\prod_{j \in J} (1-w/t_j)   }{\prod_{j \in I} (1-w/t_j) }
 {\rm d}w 
  \quad (|w/t_j| <1).
\end{align*}
Here we have used the Taylor expansion,
$\log (1-u)= -\sum_{k=1}^\infty u^k/k$
valid for $|u|<1$. 
Suppose $z w=1$.
Then we observe that
\[
\Omega_2=
z^{-n} \frac{\prod_{j \in J} (1-1/t_j z)   }{\prod_{j \in I} (1-1/t_j z) }
\left( - \frac{{\rm d}z}{z^2} \right)
=
- z^{|I|-|J|-m-n-2}\frac{\prod_{j \in I}(-t_j)}{\prod_{j \in J} (-t_j)}
\Omega_1.
\]
Consequently, both integrands of 
(\ref{eq:mn_uc_1}) and (\ref{eq:mn_uc_2})
coincide up to constant multiplication
if the condition
$|I|-|J|=m+n+2$
is fulfilled.
In this case the integrand of 
(\ref{eq:mn_uc_1})
reads
\[
F(z)=z^m\frac{\prod_{j \in J} (1-t_j z)}{ \prod_{j \in I} (1-t_j z)} 
\tau_{0,0} ({\boldsymbol x'}+[z^{-1}],{\boldsymbol y'}+[z]) 
\tau_{m,n}({\boldsymbol x}-[z^{-1}],{\boldsymbol y}-[z]).
\]
Since $\tau_{0,0}({\boldsymbol x},{\boldsymbol y})$ and  $\tau_{m,n}({\boldsymbol x},{\boldsymbol y})$ are entire,
$F(z)$ has the $|I|+2$ singularities:
$z=1/t_i$ (simple poles) for $i \in I$
and $z=0,\infty$ (which may be essential singularities).
Hence (\ref{subeq:mn_uc}) becomes
\begin{equation} \label{eq:integ}
\int_{C_1} F(z) {\rm d}z=\int_{C_2} F(z) {\rm d}z=0,
\end{equation}
where the integration contour $C_1$ (resp. $C_2$) is
a positively oriented small circle around $z=0$ (resp. $z=\infty$)
such that all the other singularities are exterior to it; see Figure~\ref{fig:cont}.
We verify through the Cauchy--Goursat theorem that
\begin{equation} \label{eq:res}
\sum_{i \in I}
\underset{z=1/t_i}{\rm{Res}}  F(z) {\rm d}z=0
\end{equation}
by canceling contribution of residues at 
$z=\infty$ and $z=0$
respectively 
to the first and second integrals in (\ref{eq:integ}).
In other words,
we have successfully avoided the residue calculus 
at possible essential singularities $z=0,\infty$ 
thanks to
 the presence of {\it two} bilinear equations
(\ref{eq:mn_uc_1}) and (\ref{eq:mn_uc_2}).

Now we prepare some notations.
For a function $f=f({\boldsymbol x},{\boldsymbol y})$, 
we define a shift operator 
$T_i$ 
by
$T_i(f)=f({\boldsymbol x}-[t_i],{\boldsymbol y}-[{t_i}^{-1}])$.
We also write 
$T_{\{i_1,i_2,\ldots,i_r \}}(f)=T_{i_1} \circ T_{i_2} \circ \cdots \circ T_{i_r}(f)$
for brevity.
Then (\ref{eq:res}) takes the following form:
\[
\sum_{i \in I} {t_i}^n \frac{ \prod_{j \in J} (t_i-t_j) }{ \prod_{j \in I \setminus \{i\}} (t_i-t_j)}
T_{I \setminus \{i\}} (\tau_{0,0}) T_{J \cup \{i\}} (\tau_{m,n}) =0,
\]
which can be regarded as a difference equation with each $t_i$ 
being the difference interval.
Along the same lines as (\ref{eq:mUC1}), 
also (\ref{eq:mUC2}) and (\ref{eq:mUC3})
generate similar difference equations.
Summarizing above we have the

\begin{prop} \label{prop:duc}
The following difference equations hold.
\begin{enumerate}
\item 
If 
$|I|-|J|=m+n+2$ and $m,n \geq 0$,
then
\begin{equation}  \label{eq:dUC1}
\sum_{i \in I} {t_i}^n \frac{ \prod_{j \in J} (t_i-t_j) }{ \prod_{j \in I \setminus \{i\}} (t_i-t_j)}
T_{I \setminus \{i\}} (\tau_{0,0}) T_{J \cup \{i\}} (\tau_{m,n}) =0.
\end{equation}

\item 
If  
$|I|-|J|=n+1$ and $n \geq 0$,
then
\begin{equation}   \label{eq:dUC2}
T_I(\tau_{0,0}) T_J(\tau_{1,n})
=
\sum_{i \in I}  
\frac{ \prod_{j \in J} (1-t_j/t_i) }{ \prod_{j \in I \setminus \{i\}} (1-t_j/t_i)}
T_{I \setminus \{i\}} (\tau_{1,0}) T_{J \cup \{i\}} (\tau_{0,n}).
\end{equation}

\item 
If
$|I|-|J|=m+1$ and $m \geq 0$, then
\begin{equation}   \label{eq:dUC3}
T_I(\tau_{0,0}) T_J(\tau_{m,1})
=
\sum_{i \in I}
\frac{ \prod_{j \in J} (1-t_i/t_j) }{ \prod_{j \in I \setminus \{i\}} (1-t_i/t_j)}
T_{I \setminus \{i\}} (\tau_{0,1}) T_{J \cup \{i\}} (\tau_{m,0}).
\end{equation}
\end{enumerate}
\end{prop}

\begin{example}\rm
Consider the case
$m=n=0$, $I=\{1,2,3\}$, and $J=\{4\}$.
Write $\tau=\tau_{0,0}$.
Then (\ref{eq:dUC1}) reduces to the equation
\begin{align*}
&(t_1-t_2)(t_3-t_4)T_{1,2}(\tau)T_{3,4}(\tau)+
(t_2-t_3)(t_1-t_4)T_{2,3}(\tau)T_{1,4}(\tau) 
\\ 
&\qquad +(t_3-t_1)(t_2-t_4)T_{1,3}(\tau)T_{2,4}(\tau)
=0,
\end{align*}
which was found by Ohta \cite{oht07}
as a quadratic relation for the universal character.

Let $m=1$, $n=0$, $I=\{1,2,3\}$, and $J=\emptyset$.
Then (\ref{eq:dUC1}) reduces to
\begin{equation}  \label{eq:dmUC}
(t_1-t_2)T_{1,2}(\tau_{0,0})T_{3}(\tau_{1,0})+
(t_2-t_3)T_{2,3}(\tau_{0,0})T_{1}(\tau_{1,0})+
(t_3-t_1)T_{1,3}(\tau_{0,0})T_{2}(\tau_{1,0})
=0.
\end{equation}
Let $n=1$, $I=\{1,2\}$, and $J=\emptyset$.
Then (\ref{eq:dUC2}) reduces to
\begin{equation} \label{eq:LdUC}
(t_1-t_2)T_{1,2}(\tau_{0,0})\tau_{1,1} 
=
t_1 T_1(\tau_{0,1}) T_2(\tau_{1,0})
-
t_2 T_2(\tau_{0,1}) T_1(\tau_{1,0}).
\end{equation}
The above difference equations, 
(\ref{eq:dmUC}) and (\ref{eq:LdUC}),
were introduced  
in a study of the connection between the universal character and
$q$-Painlev\'e equations;
see \cite{tsu05b, tsu09a}.
\end{example}

Furthermore, we can obtain
a functional  equation that involves
derivative terms
from the difference equations
through a limit process causing a confluence of the poles $z=1/t_i$.
For instance let us take the limit $t_3 \to t_1$ in (\ref{eq:dmUC}).
Rewrite $(t,s)=(t_1,t_2)$ and shift the variables as
${\boldsymbol x}\mapsto {\boldsymbol x}+[t]$
and 
${\boldsymbol y}\mapsto {\boldsymbol y}+[t^{-1}]$.
Then we find 
\begin{align} \label{eq:ddmUC_1}
& \left(
D_{\delta_t} +\frac{t}{s-t}
\right)
\tau_{0,0} ({\boldsymbol x}-[s],  {\boldsymbol y}-[s^{-1}])  
\cdot 
\tau_{1,0} ({\boldsymbol x},{\boldsymbol y} ) 
\\
& \qquad
+\frac{t }{t-s}
\tau_{0,0} ({\boldsymbol x}-[t],{\boldsymbol y}-[t^{-1}])  
\tau_{1,0} ({\boldsymbol x}+[t]-[s], {\boldsymbol y}+[t^{-1}]-[s^{-1}]) 
=0.
\nonumber
\end{align}
Here we have introduced the vector fields
\[
\delta_t=\sum_{n=1}^\infty
\left( t^n \frac{\partial}{\partial x_n}
- t^{-n} \frac{\partial}{\partial y_n}\right)
\quad
\text{and} 
\quad
\tilde{\delta}_t=\sum_{n=1}^\infty
\left(n t^n \frac{\partial}{\partial x_n}
+n  t^{-n} \frac{\partial}{\partial y_n}\right),
\]
and let $D_{\boldsymbol v}$ denote
the Hirota differential 
with respect to a vector field ${\boldsymbol v}$.
If we take continuously the limit $s \to t$ in 
(\ref{eq:ddmUC_1}) with divided by $t-s$, 
then we obtain 
\begin{equation}
\label{eq:ddmUC_2}
\left(
{D_{\delta_t}}^2 -D_{\delta_t}+D_{\tilde{\delta}_t}
\right)
\tau_{0,0} ({\boldsymbol x}-[t],  {\boldsymbol y}-[t^{-1}])  
\cdot 
\tau_{1,0} ({\boldsymbol x},{\boldsymbol y} ) 
=0.
\end{equation}

In this manner we can produce various functional equations from the UC hierarchy.
We list the ones 
relevant to the following sections.

\begin{prop}  \label{prop:muc}
The following difference {\rm(}and differential{\rm)}
equations hold{\rm:}
\begin{align}  
 \label{eq:muc_a}
 &(t-s) \tau_{m,n} ({\boldsymbol x}-[t]-[s],{\boldsymbol y}-[t^{-1}]-[s^{-1}])  
\tau_{m+1,n+1} ({\boldsymbol x},{\boldsymbol y}) \\
\nonumber
& \qquad 
- t \tau_{m,n+1} ({\boldsymbol x}-[t], {\boldsymbol y}-[t^{-1}])  
\tau_{m+1,n} ({\boldsymbol x}-[s],{\boldsymbol y}-[s^{-1}])  \\
& \qquad
+s \tau_{m,n+1} ({\boldsymbol x}-[s],{\boldsymbol y}-[s^{-1}]) 
\tau_{m+1,n} ({\boldsymbol x}-[t],{\boldsymbol y}-[t^{-1}])=0,
\nonumber
 \\
 \label{eq:muc_b}
 &\left( D_{\delta_t}-1 \right) \tau_{m,n+1}({\boldsymbol x},{\boldsymbol y}) 
 \cdot 
 \tau_{m+1,n}({\boldsymbol x},{\boldsymbol y}) 
 \\
\nonumber
& \qquad + \tau_{m,n} ({\boldsymbol x} -[t] , {\boldsymbol y} -[t^{-1}] ) 
\tau_{m+1,n+1} ({\boldsymbol x} +[t],{\boldsymbol y} +[t^{-1}] )
=0,
 \\
 &   \label{eq:muc_c}
 \left(
D_{\delta_t} +\frac{t}{s-t}
\right)
\tau_{m,n} ({\boldsymbol x}-[s],  {\boldsymbol y}-[s^{-1}])  
\cdot 
\tau_{m+1,n} ({\boldsymbol x},{\boldsymbol y} ) 
\\
& \qquad
+\frac{t }{t-s}
\tau_{m,n} ({\boldsymbol x}-[t],{\boldsymbol y}-[t^{-1}])  
\tau_{m+1,n} ({\boldsymbol x}+[t]-[s], {\boldsymbol y}+[t^{-1}]-[s^{-1}]) 
=0.
\nonumber
\end{align}
\end{prop}

\pf
Clearly
 (\ref{eq:muc_a}) and (\ref{eq:muc_c})
 are equivalent to  (\ref{eq:LdUC})
 and (\ref{eq:ddmUC_1}),
 respectively.
Taking the limit $s \to t$ in  (\ref{eq:muc_a}) 
leads to (\ref{eq:muc_b}). 
\qed

\section{Homogeneous $\tau$-sequence and its Weyl group symmetry}
\label{sect:hom}

This section is concerned with a sequence of homogeneous solutions of the UC hierarchy,
connected by vertex operators.
We show that such a sequence naturally admits a commutative pair of Weyl group actions of type $A$
generated by a permutation of two serial vertex operators.

We first introduce partial differential operators
$V_X(c)$ and  $V_Y(c')$ ($c, c' \in {\mathbb C}$)
defined by 
\[
V_X(c)= \int_\gamma X^+(z)z^{-c-1} {\rm d}z
\quad \text{and} \quad
V_Y(c')= \int_{\gamma'} Y^+(z^{-1})z^{c'-1} {\rm d}z,
\]
where the integration paths 
$\gamma, \gamma':[0,1]\to {\mathbb C}$
is taken such that $[ X^+(z)z^{-c}]_{\gamma(0)}^{\gamma(1)}
=[ Y^+(z^{-1})z^{c'}]_{\gamma'(0)}^{\gamma'(1)}=0$.
For instance $\gamma$ and $\gamma'$
can be chosen to be cycles.
Note that in general
$\gamma$ and $\gamma'$ may depend on $c$ and $c'$, respectively.
It is easy to see that
$V_X(c)$ and  $V_Y(c')$ mutually commute.

Suppose $\tau_{0,0}=\tau_{0,0}({\boldsymbol x},{\boldsymbol y})$
to be a solution of the UC hierarchy 
(\ref{eq:UCH})
satisfying the homogeneity
$E\tau_{0,0}=d_{0,0}\tau$.
Instead of (\ref{eq:seq}),
let us consider a sequence
$\{ \tau_{m,n} \}_{m,n \geq 0}$
determined recursively by
\[
\tau_{m+1,n}=V_X(c_m)\tau_{m,n} \quad \text{and}
\quad
\tau_{m,n+1}=V_Y(c'_n)\tau_{m,n}
\]
for arbitrary constant parameters $c_m,c'_n \in{\mathbb C}$ given.
Since the UC hierarchy (\ref{eq:UCH})
takes the form of bilinear equations,
it can be verified in exactly the same way as (\ref{eq:seq})
that each $\tau_{m,n}$ gives a solution of (\ref{eq:UCH}).
Furthermore, they all obey the homogeneity
\[
E \tau_{m,n}=d_{m,n} \tau_{m,n}
\]
with
$d_{m+1,n}=d_{m,n}+c_m$ and $d_{m,n+1}=d_{m,n}-c'_n$,
as a consequence of the formulae 
$[E,V_X(c)]=c V_X(c)$ and
$[E,V_Y(c')]=-c' V_Y(c')$; cf. \cite[Lemma~2.4]{tsu09b}.
Hence the balancing condition
\[
d_{m,n}+d_{m+1,n+1}=d_{m,n+1}+d_{m+1,n}
\]
is fulfilled.
We call 
the above sequence of homogeneous solutions of the UC hierarchy
a {\it homogeneous $\tau$-sequence}.
Obviously, any functional equation in Sect.~\ref{sect:fe} still remains valid for the homogeneous $\tau$-sequence $\{\tau_{m,n}\}$;
we may call also $V_X(c)$ and  $V_Y(c')$ 
vertex operators.

\begin{example}  \label{ex:uc}
\rm
If we take $c=c'=n$ to be an integer and each $\gamma$ and $\gamma'$ a positively oriented small circle around the origin $z=0$,
then $V_X(n)= 2 \pi \sqrt{-1}X_n^+$ and $V_Y(n)= 2 \pi \sqrt{-1} Y_n^+$ according to 
(\ref{subeq:vo}).
Recall now that these operators
play roles of raising operators for the universal characters;
namely, 
\[
S_{[\lambda,\mu]}({\boldsymbol x},{\boldsymbol y})=
X^+_{\lambda_1} \ldots X^+_{\lambda_\ell} Y^+_{\mu_1} \ldots Y^+_{\mu_{\ell'}} . 1
\]
for any pair of partitions
$\lambda=(\lambda_1,\lambda_2,\ldots, \lambda_\ell)$
and
$\mu=(\mu_1,\mu_2,\ldots,\mu_{\ell'})$;
see \cite[Theorem~1.2]{tsu04}.
Starting from a trivial solution $\tau({\boldsymbol x},{\boldsymbol y})  =S_{[\emptyset,\emptyset]}({\boldsymbol x},{\boldsymbol y}) \equiv 1$ of the UC hierarchy,
we thus obtain a homogeneous $\tau$-sequence expressed in
the universal characters by
successive application of $X_n^+$ and $Y_n^n$.
\end{example}

Next we consider the Weyl group symmetry of the homogeneous 
$\tau$-sequence.
Fix a positive integer $k$.
Let us look at $m=k$ sites
in the $(m,n)$-lattice
and interchange 
the $(k-1)$th and $k$th operations of the vertex operators $V_X$
in view of the fermionic relation 
$V_X(a)V_X(b)+V_X(b-1)V_X(a+1)=0$.
To be more precise,
we transform the original sequence
\[
\cdots
\stackrel{V_X(c_{k-2})}{\longrightarrow} 
\tau_{k-1,n} \stackrel{V_X(c_{k-1})}{\longrightarrow} \tau_{k,n} \stackrel{V_X(c_{k})}{\longrightarrow} \tau_{k+1,n} \stackrel{V_X(c_{k+1})}{\longrightarrow} \cdots
\]
into a new one
\[
\cdots  
\stackrel{V_X(c_{k-2})}{\longrightarrow} 
\tau_{k-1,n} \stackrel{V_X(c_{k}+1)}{\longrightarrow} \hat{\tau}_{k,n} 
\stackrel{V_X(c_{k-1}-1)}{\longrightarrow} \tau_{k+1,n} \stackrel{V_X(c_{k+1})}{\longrightarrow} \cdots
\]
that is identical with the original one 
except 
$\tau_{k,n}$ is replaced by 
\[
\hat{\tau}_{k,n} =V_X(c_k+1) \tau_{k-1,n}.
\]
Besides,
the degree of $\hat{\tau}_{k,n}$ reads 
\[
\hat{d}_{k,n}=d_{k-1,n}+c_{k}+1
=d_{k-1,n}-d_{k,n}+d_{k+1,n}+1.
\]
We refer to the above permutation of vertex operators 
as $r_k$.
Put 
\begin{equation} \label{eq:alpha}
\alpha_{k}=  \hat{d}_{k,n}-d_{k,n}=d_{k-1,n}- 2 d_{k,n}+d_{k+1,n} +1,
\end{equation}
which is a quantity that does not depend on $n$.
The operation $r_k$ 
induces 
the transformation
\[
r_k(\alpha_{k})= -\alpha_{k}, \quad 
r_k(\alpha_{k\pm 1})= \alpha_{k \pm 1}+\alpha_{k},
\quad  \text{and} \quad
r_k(\alpha_{\ell})= \alpha_{\ell} \quad (\ell \neq k, k \pm 1).
\]
Therefore $\alpha_{k}$ can be regarded as a root variable
of  the Weyl group of type $A$,
and
$\langle r_k \rangle$ indeed fulfills 
its fundamental relations
\[
{r_k}^2=1, \quad 
r_k r_{k \pm 1} r_k = r_{k \pm  1} r_k  r_{k \pm 1},
\quad
\text{and} \quad
r_k r_\ell=r_\ell r_k  \quad (\ell \neq k, k \pm 1).
\]
Along the same lines
we can derive from a permutation of operators
$V_Y$
another action of the Weyl group of type $A$,
which commutes with the previous one.
As demonstrated in Sect~\ref{sect:weyl},
this kind of Weyl group actions gives rise to 
a group of
birational canonical transformations
of the Hamiltonian system ${\cal H}_{L,N}$.

We  conclude this section with some formulae
that will be employed later.

\begin{lemma}
It holds that
\begin{align}  \label{eq:bil_hat_1}
&\tau_{k-1,n} \otimes \tau_{k+1,n}- \sum_{i+j=-1} X_i^- \hat{\tau}_{k,n} \otimes X_j^+ \tau_{k,n}
= 
\sum_{i+j=-1} Y_i^- \hat{\tau}_{k,n} \otimes Y_j^+ \tau_{k,n}
=0,
\\
\label{eq:bil_hat_2}
&\tau_{k-1,n+1} \otimes \tau_{k+1,n}- \sum_{i+j=-1} X_i^- \hat{\tau}_{k,n+1} \otimes X_j^+ \tau_{k,n}
=\hat{\tau}_{k,n} \otimes \tau_{k,n+1}- 
\sum_{i+j=0} Y_i^- \hat{\tau}_{k,n+1} \otimes Y_j^+ \tau_{k,n}
=0.
\end{align}
\end{lemma}

\pf 
First we have
\begin{equation} \label{eq:10_muc}
\sum_{i+j=-2} X_i^- \tau_{k-1,n} \otimes X_j^+ \tau_{k,n}
=
\sum_{i+j=-1} Y_i^- 
\tau_{k-1,n} \otimes Y_j^+ \tau_{k,n}
=0,
\end{equation}
which is equivalent to 
(\ref{eq:mUC1}) with $m=1$ and $n=0$.
Applying $V_X(c_k+2)\otimes 1$ 
and  $V_X(c_k+1)\otimes 1$
respectively to the first and second equalities in (\ref{eq:10_muc})
leads to (\ref{eq:bil_hat_1}). 
We deduce (\ref{eq:bil_hat_2}) from (\ref{eq:bil_hat_1}) by applying 
$V_Y(c'_n)\otimes 1$.
\qed

\begin{lemma}
The following difference {\rm(}and differential{\rm)}
equations hold{\rm:}
\begin{align}
 &  \label{eq:diff_hat_1}
 D_{\delta_t}  \hat{\tau}_{k,n}({\boldsymbol x},{\boldsymbol y})
 \cdot
 \tau_{k,n}({\boldsymbol x},{\boldsymbol y})
 -
 t \tau_{k-1,n}({\boldsymbol x}-[t],{\boldsymbol y}-[t^{-1}]) 
\tau_{k+1,n}({\boldsymbol x}+[t],{\boldsymbol y}+[t^{-1}])=0,
\\
&  \label{eq:diff_hat_2}
t \tau_{k-1,n+1}({\boldsymbol x}-[t],{\boldsymbol y}-[t^{-1}]) 
\tau_{k+1,n}({\boldsymbol x},{\boldsymbol y})
\\
& \qquad
- \hat{\tau}_{k,n+1}({\boldsymbol x},{\boldsymbol y})
 \tau_{k,n}({\boldsymbol x}-[t],{\boldsymbol y}-[t^{-1}])
  +\hat{\tau}_{k,n}({\boldsymbol x}-[t],{\boldsymbol y}-[t^{-1}])
  \tau_{k,n+1}({\boldsymbol x},{\boldsymbol y})=0.
  \nonumber
\end{align}
\end{lemma}

\pf
The verification can be done along the same argument as 
Proposition~\ref{prop:duc}.
First we shall regard the symbol $f \otimes g$ as a product 
of two functions 
$f({\boldsymbol x'},{\boldsymbol y'})g({\boldsymbol x},{\boldsymbol y})$
in distinct indeterminates
$({\boldsymbol x'},{\boldsymbol y'})$ and $({\boldsymbol x},{\boldsymbol y})$.
Taking the variables in (\ref{eq:bil_hat_1})
as ${\boldsymbol x}-{\boldsymbol x'}=2[t]$ and ${\boldsymbol y}-{\boldsymbol y'}=2[t^{-1}]$
thus leads to (\ref{eq:diff_hat_1}).
Similarly, we deduce (\ref{eq:diff_hat_2}) from 
(\ref{eq:bil_hat_2})
with
${\boldsymbol x}-{\boldsymbol x'}=[t]$ and ${\boldsymbol y}-{\boldsymbol y'}=[t^{-1}]$.
\qed

\section{Similarity reduction of UC hierarchy}
\label{sect:sim}

In this section
we consider
a reduction of the UC hierarchy
by requiring certain homogeneity and periodicity. 
As a result we derive a finite-dimensional integrable system of
partial differential equations,
denoted by ${\cal G}_{L,N}$,
which provides
an extension of both the Garnier system 
and the sixth Painlev\'e equation $P_{\rm VI}$.

Fix integers $L \geq 2$ and $N \geq 1$. 
Let $\tau_{m,n}=\tau_{m,n}({\boldsymbol x},{\boldsymbol y})$ 
be a sequence of the solutions satisfying
(\ref{eq:muc_a})--(\ref{eq:muc_c}) in Proposition~\ref{prop:muc}.
Suppose that
$\tau_{m,n}$ are homogeneous of degree $d_{m,n} \in {\mathbb C}$,
i.e.,
\[E \tau_{m,n}=d_{m,n} \tau_{m,n}
\quad \text{with} \quad
E=\sum_{n=1}^\infty 
\left(n x_n \frac{\partial}{\partial x_n} - n y_n \frac{\partial}{\partial y_n}
\right),
\]
and fullfil
the periodic condition:
$\tau_{m+L,n}=\tau_{m,n+L}=\tau_{m,n}$
(up to multiplication by constants).
Remark here that the relation
$d_{m,n}+d_{m+1,n+1}=d_{m,n+1}+d_{m+1,n}$
necessarily holds; 
cf. Sect.~\ref{sect:hom}.
Let us replace the independent variables $x_n$ and $y_n$ 
respectively with
the $n$th and $(-n)$th power sum
of new ones ${\boldsymbol t}=(t_0,t_1,\ldots,t_N)$ as
\begin{equation}  \label{eq:special_gar}
x_n=\frac{1}{n} \sum_{i=0}^N \theta_i {t_i}^n
 \quad \text{and} \quad 
y_n=\frac{1}{n}\sum_{i=0}^N \theta_i {t_i}^{-n}.  
\end{equation}
Consequently we have
\begin{align}  
t_i \frac{\partial}{\partial t_i}&=
t_i \sum_{n=1}^\infty  
 \left(
 \frac{\partial x_n}{\partial t_i} \frac{\partial}{\partial x_n}+
 \frac{\partial y_n}{\partial t_i} \frac{\partial}{\partial y_n}
 \right)
 =\theta_i \sum_{n=1}^\infty  \left( {t_i}^n \frac{\partial}{\partial x_n}- {t_i}^{-n} \frac{\partial}{\partial y_n}\right)
 =\theta_i \delta_{t_i},
\label{eq:tdt}
 \\
E&=\sum_{i=0}^N t_i \frac{\partial}{\partial t_i}
=  \sum_{i=0}^N \theta_i  \delta_{t_i}.
\label{eq:E=}
\end{align}
In view of the homogeneity,
no generality is lost by taking
$t_0=1$.
Set 
$\sigma_{m,n}({\boldsymbol \theta},{\boldsymbol t})=\tau_{m,n}({\boldsymbol x},{\boldsymbol y})$
under the above conditions.
For notational simplicity we shall use the abbreviation $\sigma_{m,n}(\theta_i \pm 1)$
to mean that 
among the constant parameters 
${\boldsymbol \theta}=(\theta_0,\theta_1,\ldots,\theta_N)$
only the indicated one
$\theta_i$ 
is shifted by $\pm 1$
while all the others are unchanged.
Then we have the

\begin{prop}  \label{prop:bil}
The functions $\sigma_{m,n}=\sigma_{m,n}({\boldsymbol \theta},{\boldsymbol t})$ satisfy the bilinear equations
\begin{subequations}  \label{subeq:bil_gar}
\begin{align}
&(t_i-t_j)\sigma_{m,n} \sigma_{m+1,n+1}(\theta_i+1,\theta_j+1)
\label{eq:bil_gar_1}
\\
&\qquad = t_i \sigma_{m+1,n} (\theta_i+1)
\sigma_{m,n+1} (\theta_j+1)
- t_j \sigma_{m+1,n} (\theta_j+1)
\sigma_{m,n+1} (\theta_i+1), 
\nonumber
\\
& \left( t_i D_{i} + \theta_i \right)
\sigma_{m+1,n} \cdot \sigma_{m,n+1}
=
\theta_i 
\sigma_{m,n} (\theta_i-1)\sigma_{m+1,n+1}(\theta_i+1),
 \label{eq:bil_gar_2}
 \\
&
\left((t_j-t_i)D_{i} +\theta_i \right) 
\sigma_{m,n}(\theta_j-1) \cdot \sigma_{m+1,n}
=
\theta_i \sigma_{m,n}(\theta_i-1)
\sigma_{m+1,n}(\theta_i+1,\theta_j-1),
 \label{eq:bil_gar_3}
\end{align}
together with the homogeneity constraint 
\begin{equation} \label{eq:hom_gar}
\sum_{i=0}^N t_i \frac{\partial   \sigma_{m,n}}{\partial t_i}
 =d_{m,n} \sigma_{m,n}.
\end{equation}
\end{subequations}
Here $D_i$ denotes the Hirota differential with respect to 
$\partial/\partial t_i$.
\end{prop}

\pf 
It is immediate to obtain
(\ref{eq:bil_gar_1}) 
from
(\ref{eq:muc_a})
with $(t,s)=(t_i,t_j)$.
Using (\ref{eq:tdt}) we verify 
 (\ref{eq:bil_gar_2}) and (\ref{eq:bil_gar_3})
 from (\ref{eq:muc_b}) and  (\ref{eq:muc_c}),
 respectively.
\qed
\ \\

Next we shall write down nonlinear differential equations
for appropriately chosen dependent variables.
Let us introduce the functions 
$f^{(i)}_{m,n}=f^{(i)}_{m,n}({\boldsymbol \theta},{\boldsymbol t})$ 
and 
$g^{(i)}_{m,n}=g^{(i)}_{m,n}({\boldsymbol \theta},{\boldsymbol t})$
defined by
\begin{align}
f^{(i)}_{m,n}
&= \frac{\sigma_{m,n-1}(\theta_i+1)  \sigma_{m-1,n-1}}{\sigma_{m-1,n}(\theta_i+1)  \sigma_{m,n-2} }, 
\label{eq:gar_f}
\\
g^{(i)}_{m,n}
&=  
\frac{t_i D_{i} \sigma_{m,n-1} \cdot \sigma_{m-1,n}}{\sigma_{m,n-1} \sigma_{m-1,n}} +\theta_i
=
\theta_i
\frac{\sigma_{m-1,n-1} (\theta_i-1) \sigma_{m,n}(\theta_i+1)
}{   \sigma_{m,n-1} \sigma_{m-1,n} }
\label{eq:gar_g}
\end{align}
for $i=0,1,\ldots,N$.
Note that the second equality in (\ref{eq:gar_g})
is a consequence of (\ref{eq:bil_gar_2}).
We have the conservation law
\begin{equation}  \label{eq:cons}
\prod_{j=1}^{L}f^{(i)}_{m+j,n-j}=1 \quad
\text{and} \quad
\sum_{j=1}^{L} g^{(i)}_{m+j,n-j}= L \theta_i.
\end{equation}
In addition we prepare
auxiliary variables 
$U_{m,n}^{(i,j)}$ and $V_{m,n}^{(i,j)}$
$(i \neq j)$
defined by
\begin{align*}
U_{m,n}^{(i,j)}
&= \frac{\theta_i t_j}{t_i-t_j} 
\frac{ \sigma_{m,n-1}(\theta_i-1,\theta_j+1)
 \sigma_{m,n}(\theta_i+1) }{ \sigma_{m,n-1} \sigma_{m,n}(\theta_j+1) },
\\
V_{m,n}^{(i,j)}
&= \frac{\theta_i t_i}{t_i-t_j} 
\frac{ \sigma_{m-1,n}(\theta_i-1,\theta_j+1)
 \sigma_{m,n}(\theta_i+1) }{ \sigma_{m-1,n} \sigma_{m,n}(\theta_j+1) }.
 \end{align*}
Then we have the following relations among the dependent variables.

 \begin{lemma}
 For $i \neq j$,
 it holds that
 \begin{align}
V_{m,n}^{(i,j)}-U_{m,n}^{(i,j)}
&=g_{m,n}^{(i)}, \label{eq:uv_fg_1}
\\
\frac{U_{m-1,n}^{(i,j)}}{ V_{m,n-1}^{(i,j)} }
&= \frac{t_j f_{m,n}^{(j)}}{t_i f_{m,n}^{(i)} },
\label{eq:uv_fg_2}
\\
V_{m,n-1}^{(i,j)}-U_{m-1,n}^{(i,j)}
&=g_{m,n}^{(i)}(\theta_j+1), 
\label{eq:uv_fg_3}
\\
\frac{U_{m,n}^{(i,j)}}{ V_{m,n}^{(i,j)} }
&= \frac{t_j f_{m,n}^{(j)} (\theta_i-1) }{t_i f_{m,n}^{(i)} (\theta_i-1) }.
\label{eq:uv_fg_4}
\end{align}
 \end{lemma}
 
 \pf
 Clearly
   (\ref{eq:uv_fg_2}) and (\ref{eq:uv_fg_4})
 are direct consequences of
 the definition of $f_{m,n}^{(i)}$, (\ref{eq:gar_f}).
We obtain
  (\ref{eq:uv_fg_1}) and (\ref{eq:uv_fg_3}) from
  the bilinear equation
  (\ref{eq:bil_gar_1}).
 \qed
 \\

Solving the linear equations 
(\ref{eq:uv_fg_1}) and (\ref{eq:uv_fg_2})
for $U_{m,n}^{(i,j)}$ and $V_{m,n}^{(i,j)}$
with the aid of the $(L,L)$-periodicity,
we conclude that 
\begin{align*}
U_{m,n}^{(i,j)}
&= \frac{1}{\left( \frac{t_i}{t_j}\right)^L -1}
\sum_{b=1}^L   g^{(i)}_{m-b+1,n+b-1}  
\prod_{a=1}^{b-1} \frac{t_i f^{(i)}_{m-a+1,n+a}}{t_j f^{(j)}_{m-a+1,n+a}}, \\
V_{m,n}^{(i,j)}
&= \frac{1}{\left( \frac{t_i}{t_j}\right)^L -1}
\sum_{b=1}^L   g^{(i)}_{m-b,n+b}  
\prod_{a=0}^{b-1} \frac{t_i f^{(i)}_{m-a,n+a+1}}{t_j f^{(j)}_{m-a,n+a+1}}.
\end{align*}
In fact
$U_{m,n}^{(i,j)}$ and $V_{m,n}^{(i,j)}$ can be expressed as  
 polynomials in $f^{(i)}_{m,n}$
and $g^{(i)}_{m,n}$
via
 (\ref{eq:cons}).

\begin{thm}   \label{thm:gar}
The functions
$f^{(i)}_{m,n}$ 
and 
$g^{(i)}_{m,n}$
satisfy the system of nonlinear differential equations
\begin{subequations} \label{subeq:gar}
\begin{align}
t_i \frac{\partial f_{m,n}^{(i)}}{\partial t_i}
&= \left(
\kappa_{m,n} -g_{m,n-1}^{(i)} 
+\sum_{j \neq i} \left( U_{m-1,n}^{(j,i)} - V_{m,n-1}^{(j,i)} \right)
\right)  f_{m,n}^{(i)},
 \label{eq:gar_a}
\\
t_j \frac{\partial f_{m,n}^{(i)}}{\partial t_j}
&=\left(-g_{m,n-1}^{(i)}-U_{m-1,n}^{(j,i)} +V_{m,n-1}^{(j,i)}
\right) f_{m,n}^{(i)}
\quad
(i \neq j),
 \label{eq:gar_b}
\\
t_i \frac{\partial g_{m,n}^{(i)}}{\partial t_i}
&= - \sum_{j \neq i}
\left(
U_{m,n}^{(i,j)}g_{m,n}^{(j)}+ V_{m,n}^{(j,i)} g_{m,n}^{(i)}
\right),
 \label{eq:gar_c}
\\
t_j \frac{\partial g_{m,n}^{(i)}}{\partial t_j}
&= U_{m,n}^{(i,j)}g_{m,n}^{(j)}+V_{m,n}^{(j,i)} g_{m,n}^{(i)}
\quad
(i \neq j),
 \label{eq:gar_d}
\end{align}
\end{subequations}
where
\begin{equation}
 \label{eq:kappa}
\kappa_{m,n}=d_{m,n-1}-d_{m-1,n}+\sum_{i=0}^N \theta_i
=\sum_{i=0}^N g_{m,n}^{(i)} \in {\mathbb C} 
\end{equation}
are constant parameters.
\end{thm}

For each $(m,n)$ fixed
the system
(\ref{subeq:gar}) is closed with respect to the
$2 L N$-tuple of dependent variables
$g^{(i)}_{m+j,n-j}$
and
$f_{m+j,n-j+1}^{(i)}/f_{m+j,n-j+1}^{(0)}$,
where  $i=1,2,\ldots,N$
and $j \in {\mathbb Z}/L {\mathbb Z}$.
Moreover, it possesses the $2N$ 
conserved quantities;
recall (\ref{eq:cons}).
Accordingly
the dimension of the phase space is essentially $2N(L-1)$.
If $L=2$ then it is in fact
equivalent to the Garnier system in $N$ variables,
whose phase space is $2N$-dimensional;
see also the appendix.
Let ${\cal G}_{L,N}$ denote the nonlinear system
(\ref{subeq:gar}).
As shown in Sect.~\ref{sect:ham}, 
the system
${\cal G}_{L,N}$ can be transformed into a canonical Hamiltonian system 
with polynomial Hamiltonian functions.
\\

\noindent
{\it Proof of Theorem~\ref{thm:gar}.}
We shall demonstrate only (\ref{eq:gar_a}) here because the others 
can be done in quite a similar manner.
By virtue of
the homogeneity (\ref{eq:hom_gar})
we see that 
\[
\sum_{i=0}^N \left(g_{m,n}^{(i)} -\theta_i \right)
=\sum_{i=0}^N \frac{ t_i D_i \sigma_{m.n-1} \cdot \sigma_{m-1,n}}{\sigma_{m.n-1} \sigma_{m-1,n}} 
=d_{m,n-1}-d_{m-1,n}.
\]
Therefore
(\ref{eq:kappa}) certainly holds.
By combining (\ref{eq:uv_fg_3}) with (\ref{eq:kappa}) we have also
\begin{align}  \label{eq:shift_g_2}
g_{m,n}^{(i)}(\theta_i+1)
&=\kappa_{m,n}+1 -  \sum_{j \neq i} g_{m,n}^{(j)}(\theta_i+1)
\\
&= \kappa_{m,n}+1+ \sum_{j \neq i} \left( U_{m-1,n}^{(j,i)} -V_{m,n-1}^{(j,i)} \right).
\nonumber 
\end{align}
Taking the logarithmic derivative of $f_{m,n}^{(i)}$
shows that
\begin{align*}
\frac{t_i}{f_{m,n}^{(i)}}\frac{\partial f_{m,n}^{(i)}}{\partial t_i}
&=
\frac{t_i D_i \sigma_{m,n-1}(\theta_i+1) \cdot \sigma_{m-1,n}(\theta_i+1) }{\sigma_{m,n-1}(\theta_i+1)  \sigma_{m-1,n}(\theta_i+1)} 
- \frac{ t_i D_i \sigma_{m,n-2} \cdot \sigma_{m-1,n-1} }{\sigma_{m,n-2} \sigma_{m-1,n-1}}
\\
&=  g_{m,n}^{(i)}(\theta_i+1) - g_{m,n-1}^{(i)} -1, 
\quad \text{using (\ref{eq:gar_g})},
\\
&=\kappa_{m,n}-g_{m,n-1}^{(i)} +
\sum_{j \neq i} \left( U_{m-1,n}^{(j,i)} -V_{m,n-1}^{(j,i)} \right),
 \quad \text{using (\ref{eq:shift_g_2})}.
\end{align*}
We have verified (\ref{eq:gar_a}) as desired. 
\qed

\begin{remark}[Toda equation] \rm
We shall derive a differential-difference 
equation of Toda-type
for 
$\sigma_{m,n}$,
associated with the shift 
$(\theta_i,\theta_j)\mapsto (\theta_i+1,\theta_j-1) $
of parameters.
First we differentiate with respect to $s$ the equation
(\ref{eq:muc_c})
after shifting the variables 
$({\boldsymbol x},{\boldsymbol y})$
to $({\boldsymbol x}+[s]/2,{\boldsymbol y}+[s^{-1}]/2)$.
We thus find that
\begin{align*}
& \nonumber
 \left(
(s-t)D_{\delta_t} D_{\delta_s} -2s D_{\delta_t} 
+t D_{\delta_s} 
\right)
\tau_{m,n} 
\left({\boldsymbol x}-\frac{[s]}{2},  {\boldsymbol y}-\frac{[s^{-1}]}{2} \right)  
\cdot 
\tau_{m+1,n}
\left({\boldsymbol x}+\frac{[s]}{2},  {\boldsymbol y}+\frac{[s^{-1}]}{2} \right) 
\\
& \qquad
+t  D_{\delta_s} 
\tau_{m,n}
\left({\boldsymbol x}-[t]+\frac{[s]}{2},  {\boldsymbol y}-[t^{-1}]+\frac{[s^{-1}]}{2} \right) 
\cdot
\tau_{m+1,n} 
\left({\boldsymbol x}+[t]-\frac{[s]}{2}, {\boldsymbol y}+[t^{-1}]-\frac{[s^{-1}]}{2} \right) 
=0.
\end{align*}
 Substitution of (\ref{eq:special_gar})
 and $(t,s)=(t_i,t_j)$  
 leads to
  \begin{align*}
& \nonumber
 \left(
(t_j-t_i) D_i D_j - (2 \theta_j+1) D_i
+\theta_i D_j
\right)
\sigma_{m,n}
\cdot
\sigma_{m+1,n}(\theta_j+1)
\\
& \qquad
+\theta_i D_j
\sigma_{m,n}(\theta_i-1,\theta_j+1)
\cdot
\sigma_{m+1,n}(\theta_i+1)=0.
 \end{align*}
 Hence, with the aid of (\ref{eq:bil_gar_3}),
 we verify that
 \begin{align}  \label{eq:toda_pre1}
\lefteqn{
(t_i-t_j)^2 \frac{ D_i D_j \sigma_{m,n}
\cdot
\sigma_{m+1,n}(\theta_j+1) 
}{\sigma_{m,n} \sigma_{m+1,n}(\theta_j+1)  }
}  \\
&\qquad
=
-\theta_i (2 \theta_j+1)
 +
\theta_i \theta_j
\frac{ 
\sigma_{m,n}(\theta_i+1,\theta_j-1) \sigma_{m,n}(\theta_i-1,\theta_j+1) }{{\sigma_{m,n}}^2} 
\nonumber \\
& \qquad \qquad +
\theta_i (\theta_j+1)
\frac{ 
\sigma_{m+1,n}(\theta_i+1) \sigma_{m+1,n}(\theta_i-1,\theta_j+2) }{{\sigma_{m+1,n}(\theta_j+1)}^2} 
\nonumber \\
&\qquad \qquad 
+(t_i-t_j)^2
\left(
\frac{ D_i \sigma_{m,n}\cdot \sigma_{m+1,n} (\theta_j+1)}{\sigma_{m,n} \sigma_{m+1,n}(\theta_j+1) }
\right)
\left(
\frac{ D_j \sigma_{m,n} \cdot \sigma_{m+1,n}(\theta_j+1) }{\sigma_{m,n}\sigma_{m+1,n}(\theta_j+1) }
\right).
\nonumber
 \end{align}
Next we express (\ref{eq:bil_gar_3}) in the form
\[
(t_i-t_j) \frac{ D_i
\sigma_{m,n} \cdot \sigma_{m+1,n}(\theta_j+1)}{\sigma_{m,n} \sigma_{m+1,n}(\theta_j+1)}
= \theta_i- \theta_i 
\frac{\sigma_{m,n}(\theta_i-1, \theta_j+1)
\sigma_{m+1,n}(\theta_i+1)}{
\sigma_{m,n} \sigma_{m+1,n}(\theta_j+1)}.
\] 
By differentiating this with respect to $t_j$,
we have
 \begin{align} \label{eq:toda_pre2}
 \lefteqn{
 (t_i-t_j)^2\frac{\partial}{\partial t_j}
 \left(
 \frac{ D_i
\sigma_{m,n} \cdot \sigma_{m+1,n}(\theta_j+1) }{\sigma_{m,n} \sigma_{m+1,n}(\theta_j+1)}
 \right)}
\\
 & \qquad
 = \theta_i 
 +
 \theta_i \theta_j
\frac{ 
\sigma_{m,n}(\theta_i+1,\theta_j-1) \sigma_{m,n}(\theta_i-1,\theta_j+1) }{{\sigma_{m,n}}^2}
\nonumber \\
& \qquad \qquad
 -
 \theta_i (\theta_j+1)
\frac{ 
\sigma_{m+1,n}(\theta_i+1) \sigma_{m+1,n}(\theta_i-1,\theta_j+2) }{{\sigma_{m+1,n}(\theta_j+1)}^2} .
 \nonumber 
 \end{align}
 Finally, combining (\ref{eq:toda_pre1}) and (\ref{eq:toda_pre2}), we arrive at the {\it Toda equation}:
 \begin{equation}  \label{eq:toda}
(t_i-t_j)^2 \frac{D_i D_j
\sigma_{m,n} \cdot \sigma_{m,n}}{  {\sigma_{m,n}}^2  }
=-2 \theta_i\theta_j +2 \theta_i\theta_j \frac{\sigma_{m,n}(\theta_i+1,\theta_j-1)\sigma_{m,n}(\theta_i-1,\theta_j+1)}{  {\sigma_{m,n}}^2  }.
\end{equation}
Note that (\ref{eq:toda}) is still valid without requiring the homogeneity and periodicity.
Such a differential-difference 
equation of Toda-type has 
previously been studied for the case of $P_{\rm VI}$,
i.e., $(L,N)=(2,1)$ 
by Okamoto \cite{oka87};
and for the case of the Garnier systems, i.e.,
$L=2$ and general $N \geq 1$,
refer to  \cite{tsu06}.
\end{remark}

\section{Particular solutions expressed in terms of the universal character}
\label{sect:ratsol}

As described in Sect.~\ref{sect:sim}, 
the system ${\cal G}_{L,N}$
is a similarity reduction of the UC hierarchy.
Since the universal characters 
$S_{[\lambda,\mu]}=S_{[\lambda,\mu]}({\boldsymbol x},{\boldsymbol y})$
are homogeneous solutions of the UC hierarchy,  
they survive through the reduction procedure;
recall Example~\ref{ex:uc}.
Therefore we can immediately construct a solution of ${\cal G}_{L,N}$
expressed in terms of the universal character.

First we recall some terminology. 
A subset ${\bf m} \subset {\mathbb Z}$ is said to be a 
{\it Maya diagram} 
if
$i \in {\bf m}$ (for $i \ll 0$) 
and
$i \notin {\bf m}$ (for $i \gg 0$).
Each Maya diagram 
${\bf m} = \{ \ldots, m_3, m_2, m_1 \}$
corresponds to
a partition
$\lambda=(\lambda_1,\lambda_2,\ldots)$ 
via $m_i-m_{i+1}=\lambda_i-\lambda_{i+1}+1$.
We can associate
with a sequence of integers
${\boldsymbol \nu}=(\nu_1,\nu_2,\ldots,\nu_{L})\in{\mathbb Z}^{L}$
a Maya diagram
\[
{\bf m}({\boldsymbol \nu})=
(L{\mathbb Z}_{<\nu_1}+1)
\cup
(L{\mathbb Z}_{<\nu_2}+2)
\cup
\cdots
\cup
(L{\mathbb Z}_{<\nu_L}+L);
\]
let $\lambda({\boldsymbol \nu})$
denote its corresponding partition.
Note that 
$\lambda ({\boldsymbol \nu}+ {\boldsymbol 1}) = \lambda({\boldsymbol \nu})$ 
where
$ {\boldsymbol 1} =(\overbrace{1,1,\ldots,1}^L) $.
We call a partition of the form $\lambda({\boldsymbol \nu})$
an {\it $L$-core} partition.
A partition $\lambda$ is 
$L$-core 
if and only if
$\lambda$ 
has no hook with length of a multiple of $L$;
see \cite[Proposition~7.13]{nou04}.
For example,
if $L=2$ and ${\boldsymbol \nu}=(0,n)$ ($n >0$) then 
the result is a staircase partition 
$\lambda({\boldsymbol \nu})=(n,n-1,\ldots,2,1)$,
thereby two-core.

There is a cyclic chain of the universal characters 
attached to $L$-core partitions
that is  connected by the action of vertex operators;
see \cite[Lemma~2.2]{tsu05a}.

\begin{lemma}
It holds that
\[
X^+_{L \nu_m -| {\boldsymbol \nu}|}
S_{[ \lambda( {\boldsymbol \nu}(m-1)), \mu]}
=\pm
S_{[ \lambda( {\boldsymbol \nu}(m)), \mu]}
\]
for arbitrary 
${\boldsymbol \nu} =(\nu_1,\nu_2,\ldots,\nu_L)
\in {\mathbb Z}^L$
and partition $\mu$.
Here
${\boldsymbol \nu}(m) ={\boldsymbol \nu} +(\overbrace{1,\ldots,1}^m ,\overbrace{0,\ldots,0}^{L-m})$
and
$|  {\boldsymbol \nu} |=\nu_1+\nu_2+\cdots+\nu_L$.
A similar formula holds for the operators $Y^+_n$ also.
\end{lemma}

Hence
we are led to the following expression of
rational solutions
of ${\cal G}_{L,N}$ 
in terms of the universal character
attached to 
a pair of $L$-core partitions.

\begin{thm}  
\label{thm:ratsol}
Let
${\boldsymbol \nu}, {\boldsymbol \nu'}  \in {\mathbb Z}^L$
be arbitrary sequences of integers.
Define 
\[
\sigma_{m,n}({\boldsymbol \theta},{\boldsymbol t})=
S_{[\lambda({\boldsymbol \nu}(m)),\lambda({\boldsymbol \nu'}(n))]}({\boldsymbol x},{\boldsymbol y}) 
\]
under the substitution
\[
x_n=\frac{1}{n} \sum_{i=0}^N \theta_i {t_i}^n
 \quad \text{and} \quad 
y_n=\frac{1}{n}\sum_{i=0}^N \theta_i {t_i}^{-n}.  
\]
Then the functions $\sigma_{m,n}$ satisfy 
the bilinear equations
{\rm(\ref{eq:bil_gar_1})}--{\rm(\ref{eq:bil_gar_3})}
and the homogeneity 
{\rm(\ref{eq:hom_gar})},
where
$d_{m,n}-d_{m-1,n}= L \nu_{m}- |{\boldsymbol \nu}|$
and
$d_{m,n}-d_{m,n-1}= -L \nu'_{n}+ |{\boldsymbol \nu'}|$.
Consequently the functions
$f_{m,n}^{(i)}$ and $g_{m,n}^{(i)}$
defined by 
{\rm(\ref{eq:gar_f})} and
{\rm (\ref{eq:gar_g})}
give a rational solution of the system ${\cal G}_{L,N}$,
{\rm (\ref{subeq:gar})}, with the parameters
$\kappa_{m,n}= L(\nu_m+\nu'_n)-|{\boldsymbol \nu}|-|{\boldsymbol \nu'}|+\sum_{i=0}^N \theta_i$. 
\end{thm}

\section{Lax formalism}
\label{sect:lax}
In this section we derive from the UC hierarchy 
the auxiliary linear problem
whose compatibility condition amounts to the nonlinear system 
 ${\cal G}_{L,N}$
 (Lax formalism). 
 It is seen that  
${\cal G}_{L,N}$ describes the
monodromy preserving deformations
of a Fuchsian system of
linear differential equations with a certain spectral type.

We introduce the {\it wave function}
\[
\psi_{m,n}({\boldsymbol x},{\boldsymbol y},k)
= \frac{ \tau_{m,n-1}({\boldsymbol x}-[k^{-1}],{\boldsymbol y}-[k])  }{ \tau_{m,n}({\boldsymbol x},{\boldsymbol y}) }e^{\xi({\boldsymbol x},k)},
\]
which is a function in 
$({\boldsymbol x},{\boldsymbol y})=(x_1,x_2,\ldots, y_1,y_2,\ldots)$
equipped with an additional parameter $k$ (the {\it spectral variable}).
Define
$\phi_{m,n}({\boldsymbol \theta}, {\boldsymbol t},k)= \psi_{m,n}({\boldsymbol x},{\boldsymbol y},k)$
under the change of variables (\ref{eq:special_gar}).
We then have the

\begin{prop}
The wave functions $\phi_{m,n}=\phi_{m,n}({\boldsymbol \theta}, {\boldsymbol t},k)$ satisfy the linear equations
\begin{align}
 \label{eq:gar_lax_1}
\phi_{m,n}
&= \frac{1}{f_{m+1,n+1}^{(i)}} \phi_{m,n+1}(\theta_i+1)
-t_i k  \phi_{m+1,n}(\theta_i+1),
\\
 \label{eq:gar_lax_2}
t_i \frac{\partial}{\partial t_i} \phi_{m,n}
&= \left(  g_{m+1,n}^{(i)}-\theta_i\right) \phi_{m,n}
+t_i k  g_{m+1,n}^{(i)}\phi_{m+1,n}(\theta_i+1),
\\
 \label{eq:gar_lax_3}
\left( k  \frac{\partial}{\partial k} - \sum_{i=0}^N t_i \frac{\partial}{\partial t_i}  \right) \phi_{m,n}
&= (d_{m,n}-d_{m,n-1}) \phi_{m,n}.
\end{align}
\end{prop}

\pf
To begin with, 
we recall
 the definition of variables 
$f_{m,n}^{(i)}$ and $g_{m,n}^{(i)}$;
see (\ref{eq:gar_f}) and (\ref{eq:gar_g}).
Substitution of $(t,s)=(t_i,1/k)$ in (\ref{eq:muc_a}) and (\ref{eq:muc_c})
produces respectively
(\ref{eq:gar_lax_1}) and (\ref{eq:gar_lax_2}),
with the aid of 
(\ref{eq:tdt}).
We deduce from the homogeneity condition
$E \tau_{m,n}=d_{m,n} \tau_{m,n}$
that 
\[
\left(E-k \frac{\partial}{\partial k}\right)\tau_{m,n}({\boldsymbol x}-[k^{-1}],{\boldsymbol y}-[k])
=d_{m,n}\tau_{m,n}({\boldsymbol x}-[k^{-1}],{\boldsymbol y}-[k]).
\]
On the other hand,  we have
$(E-k\partial/\partial k)e^{\xi({\boldsymbol x},k)}=0$.
Hence we are led to the formula
\[
\left(E- k\frac{\partial}{\partial k}\right)
 \psi_{m,n}=(d_{m,n-1}-d_{m,n})\psi_{m,n},
 \]
which implies (\ref{eq:gar_lax_3})
via (\ref{eq:E=}).
The proof is now complete.
\qed
\ \\

Because of the $(L,L)$-periodicity,
the linear equations (\ref{eq:gar_lax_1}) 
can be solved for $\phi_{m,n}(\theta_i+1)$;
thus,
\[
\phi_{m,n}(\theta_i+1)= 
\frac{1}{1-(t_i k)^L}
\sum_{b=1}^{L} (t_i k)^{b-1}  
\left(
 \prod_{a=1}^b  f_{m+a,n-a+1}^{(i)}
\right)
\phi_{m+b-1,n-b} .
\]
If we eliminate $\phi_{m+1,n}(\theta_i+1)$
from (\ref{eq:gar_lax_2}) by using
the above formula, 
then we have 
\begin{equation} \label{eq:gar_lax_2'}
t_i \frac{\partial}{\partial t_i} \phi_{m,n}
= \left(  g_{m+1,n}^{(i)}-\theta_i\right) \phi_{m,n}
+ \frac{g_{m+1,n}^{(i)}}{1-(t_i k)^L}
\sum_{b=1}^{L} (t_i k)^{b}  
\left(
 \prod_{a=1}^b  f_{m+a+1,n-a+1}^{(i)}
\right)
\phi_{m+b,n-b} .
\end{equation}
Notice that for each $m$ and $n$ fixed (\ref{eq:gar_lax_2'}) is closed with respect to $\phi_{m+j,n-j}$ $(j \in {\mathbb Z}/ L{\mathbb Z})$.
With this fact in mind,
we shall write down the linear differential equations satisfied by the vector
\[\Phi={}^{\rm T}\left( \phi_{-1,0}, k \phi_{0,-1},   k^2 \phi_{1,-2}, \ldots, k^{L-1}  \phi_{L-2, -L +1} \right).
\]
Consider the change of variables
\begin{equation}
\label{eq:zu}
z=k^L \quad
\text{and}  \quad
u_i={t_i}^{-L}.
\end{equation}
We can express
(\ref{eq:gar_lax_2'}) 
in the $L \times L$ matrix equation
\begin{equation}  \label{eq:gar_lax_B}
\frac{\partial}{\partial u_i} \Phi
=  B_i \Phi
\end{equation}
with
\begin{align*}
B_i&= {\rm diag} \left( \frac{\theta_i}{L u_i} -v_{n,n}^{(i)} \right)_{0 \leq n \leq L-1}
\\
& \quad
+ \frac{1}{z-u_i}
 \left(  
   \begin{array}{cccc}
   0 & v_{0,1}^{(i)}   &    \cdots & v_{0,L-1}^{(i)}  
   \\
       & 0&  \ddots    & \vdots
   \\      
          &        & \ddots & v_{L-2,L-1}^{(i)} 
   \\
          &        &          & 0
   \end{array}
 \right) 
 +\frac{z}{z-u_i} 
  \left(  
   \begin{array}{cccc}
   v_{0,0}^{(i)}  &    & &  \text{\LARGE $O$} 
   \\
    v_{1,0}^{(i)}    &v_{1,1}^{(i)}  &        & 
   \\      
     \vdots     &        & \ddots& 
   \\
  v_{L-1,0}^{(i)}   & v_{L-1,1}^{(i)} &\cdots& v_{L-1,L-1}^{(i)} 
   \end{array}
 \right),
\end{align*}
where
\begin{equation} \label{eq:v}
 v_{n,n+b}^{(i)}=
\frac{g_{n,-n}^{(i)}}{L}
 \prod_{a=1}^b t_i f_{n+a,-n-a+1}^{(i)}
\end{equation}
for $0 \leq n \leq L-1$ and $1 \leq b \leq L$.
Remark that
the suffix of each variable should be suitably regarded as an element of
${\mathbb Z}/L{\mathbb Z}$.

Similarly, we obtain from (\ref{eq:gar_lax_3}) 
the linear differential equation with respect to $z$:
\begin{equation}  \label{eq:gar_lax_A}
\frac{\partial \Phi}{\partial z}=A \Phi
=
\sum_{i=0}^{N+1 }\frac{A_i}{z-u_i}  \Phi,
\end{equation}  
where
$u_{N+1}=0$ and
the $L \times L$ matrices $A_i$
read
\begin{align*}
 A_i &= 
- \left(  
   \begin{array}{cccc}
   0 & v_{0,1}^{(i)}   &    \cdots & v_{0,L-1}^{(i)}  
   \\
       & 0&  \ddots    & \vdots
   \\      
          &        & \ddots & v_{L-2,L-1}^{(i)} 
   \\
          &        &          & 0
   \end{array}
 \right) 
 - u_i 
  \left(  
   \begin{array}{cccc}
   v_{0,0}^{(i)}  &    & &  \text{\LARGE $O$} 
   \\
    v_{1,0}^{(i)}    &v_{1,1}^{(i)}  &        & 
   \\      
     \vdots     &        & \ddots& 
   \\
  v_{L-1,0}^{(i)}   & v_{L-1,1}^{(i)} &\cdots& v_{L-1,L-1}^{(i)} 
   \end{array}
 \right)
 \quad (0 \leq i \leq N),
\\
A_{N+1} &= 
 \left(  
   \begin{array}{cccc}
   e_0 &   w_{0,1}      & \cdots  & w_{0,L-1}
   \\
        & e_1 &  \ddots       & \vdots
   \\     
          &        &\ddots& w_{L-2,L-1}
   \\
          &        &          & e_{L-1}    
   \end{array}
 \right)
\end{align*}
with
\[
e_n= \frac{d_{n,-n-1}-d_{n-1,-n-1}+n}{L}
\quad \text{and} \quad
w_{m,n}=  \sum_{i=0}^N v_{m,n}^{(i)}.
\]
The linear differential equation (\ref{eq:gar_lax_A}) is Fuchsian 
and has the $N+3$ regular singularities
$u_0,u_1,\ldots ,u_N$, 
$u_{N+1}=0$, 
$u_{N+2}=\infty$.
Observe that
every $A_{i}$ $(0 \leq i \leq N)$
is not full rank unlike 
$A_{N+1}$ and $A_{N+2}=-\sum_{i=0}^{N+1} A_i$.
To be specific, 
if we prepare the column vector
${\boldsymbol b}^{(i)}$ and the row vector ${\boldsymbol c}^{(i)}$
defined by
\begin{align*}
{}^{\rm T}{\boldsymbol b}^{(i)}
&=
\left( \frac{ -g_{n,-n}^{(i)} }{  L {t_i}^n  \prod_{m=1}^n  f_{ m,-m+1 }^{(i)} }
\right)_{0 \leq n \leq L-1}
=
\frac{-1}{L}
\left(
g_{0,0}^{(i)}, \frac{ g_{1,-1}^{(i)} }{ t_i f_{1,0}^{(i)}},
 \frac{ g_{2,-2}^{(i)} }{ {t_i}^2 f_{1,0}^{(i)} f_{2,-1}^{(i)}},
 \ldots,
  \frac{ g_{L-1,-L+1}^{(i)} }{ {t_i}^{L-1} f_{1,0}^{(i)} f_{2,-1}^{(i)}  \cdots f_{L-1,2}^{(i)} }
\right),
\\
{\boldsymbol c}^{(i)}
&=\left( {t_i}^n \prod_{m=1}^n  f_{ m,-m+1 }^{(i)}
\right)_{0 \leq n \leq L-1}
=
\left(1,t_i f_{1,0}^{(i)}, {t_i}^2  f_{1,0}^{(i)}f_{2,-1}^{(i)}, \ldots
, {t_i}^{L-1}  f_{1,0}^{(i)}f_{2,-1}^{(i)} \cdots f_{L-1,2}^{(i)} 
\right),
\end{align*}
then we have indeed
\begin{equation} \label{eq:A=bc}
A_{i}=
{\boldsymbol b}^{(i)} \cdot {\boldsymbol c}^{(i)} \quad \text{and} \quad
{\boldsymbol c}^{(i)} \cdot {\boldsymbol b}^{(i)}=
- \sum_{n=0}^{L-1} \frac{g^{(i)}_{n,-n}}{L}
 =-\theta_i \in {\mathbb C}
\end{equation}
for $0 \leq i \leq N$.
The matrix
$A_{N+2}=-\sum_{i=0}^{N+1} A_i$
is lower triangular and its diagonal entries are
\[
\sum_{i=0}^{N}u_i v_{n,n}^{(i)}-e_n 
=\sum_{i=0}^N \frac{ g_{n,-n}^{(i)} }{L} -e_n
=\kappa_n-e_n
\]
for $0 \leq n \leq L-1$.
Here we have used
 (\ref{eq:zu}) and (\ref{eq:v})
and put 
\begin{equation}  \label{eq:kappa_n}
\kappa_n = \frac{\kappa_{n,-n}}{L}
=\frac{d_{n,-n-1}-d_{n-1,-n}+\sum_{i=0}^{N} \theta_i}{L}
=\sum_{i=0}^N \frac{ g_{n,-n}^{(i)} }{L};
\end{equation}
 cf. (\ref{eq:kappa}).
Hence
the characteristic exponents of (\ref{eq:gar_lax_A})
at each
singularity
$z=u_i$,
i.e., the eigenvalues of each residue matrix $A_i$,
are listed 
in the following table (Riemann scheme):
\begin{equation}
\label{eq:table}
\begin{array}{cc}
\hline
\text{Singularity} & \text{Exponents}  \\ \hline
u_i  \mbox{\ } (0 \leq i \leq N )& (-\theta_i, 0, \ldots,0)  \\
u_{N+1}=0 & (e_0,e_1,\ldots,e_{L-1})  \\     
u_{N+2}=\infty & (\kappa_{0}-e_0,\kappa_{1}-e_1, \ldots,\kappa_{L-1}-e_{L-1} )  \\ 
\hline
   \end{array}
\end{equation}
Note that the relations
\begin{equation} \label{eq:rel_exp}
\sum_{n=0}^{L-1}e_n=\frac{L-1}{2}
\quad \text{and} \quad
\sum_{n=0}^{L-1}\kappa_n
=\sum_{i=0}^N \theta_i
\end{equation}
hold among the exponents.
The sum of all the exponents certainly equals zero
(Fuchs relation).

Compatibility between the above two linear equations,
(\ref{eq:gar_lax_B}) and (\ref{eq:gar_lax_A}),
is {\it a priori} established 
because both originate from 
the same bilinear equation
(\ref{eq:UCH}).
The former, (\ref{eq:gar_lax_B}), governs 
the monodromy preserving
deformation of the latter,
 (\ref{eq:gar_lax_A}), along a deformation parameter $u_i$.
The nonlinear system ${\cal G}_{L,N}$,
(\ref{subeq:gar}), 
can be recovered from the integrability condition
$\left[ \frac{\partial}{\partial u_i}-B_i ,\frac{\partial}{\partial z}-A\right]=0$
of the linear system  (\ref{eq:gar_lax_B}) and (\ref{eq:gar_lax_A}).

\begin{remark} \label{remark:scf}
\rm
In general,
we can associate with
an
$L \times L$ Fuchsian system
\begin{equation} \label{eq:scf}
\frac{\partial \Phi}{\partial z}
= A \Phi
=\sum_{i=0}^{N+1}\frac{A_{i}}{z-u_i} \Phi
\end{equation}
having $N+3$ regular singularities
$u_0,u_1,\ldots,u_N, u_{N+1}=0, u_{N+2}=\infty$
an $(N+3)$-tuple 
\[
{\cal M}=
\{
(\mu_{0,1},\mu_{0,2},\ldots,\mu_{0,\ell_0}),
(\mu_{1,1},\mu_{1,2},\ldots,\mu_{1,\ell_1}),
\ldots,
(\mu_{N+2,1},\mu_{N+2,2},\ldots,\mu_{N+2,\ell_{N+2}})
\}
\]
of partitions of $L$,
called the {\it spectral type},
in such a way that
each residue matrix $A_{i}$ has 
the eigenvalues of multiplicity $\mu_{i,j}$.
The number of {\it accessory parameters}, 
i.e., coordinates of the space of Fuchsian systems (\ref{eq:scf}) with 
given data of eigenvalues of $A_i$,
is known to be an even
\[
(N+1)L^2- \sum_{i=0}^{N+2} \sum_{j=1}^{\ell_i} {\mu_{i,j}}^2 +2.
\]

We turn now to our case.
The spectral type of 
(\ref{eq:gar_lax_A})
reads
\begin{equation}  \label{eq:spec}
 \underbrace{(L-1,1), \ldots, (L-1,1)}_{N+1},
(1,1,\ldots,1),(1,1,\ldots,1)
\end{equation}
according to its Riemann scheme (\ref{eq:table}).
Applying the above formula
we find the number of accessory parameters to be $2N(L-1)$,
which certainly equals the essential dimension of the phase space of
${\cal G}_{L,N}$
as was calculated in Sect.~\ref{sect:sim}.
\end{remark}

\begin{remark} \rm
Thanks to the algorithm proposed by Oshima \cite{osh08},
Fuchsian systems of the form (\ref{eq:scf})
with a fixed number $p$ 
of accessory parameters
can be classified by the spectral types.
Let us here take our interest 
in the Fuchsian systems that 
have four or more singularities
because they admit the monodromy preserving deformations.
If $p=2$ then 
we have a single fundamental system
whose spectral type is
$\{(1,1)^4\}=\{(1,1),(1,1),(1,1),(1,1)\}$;
and
 its deformation equation turns out to be $P_{\rm VI}$ $(={\cal G}_{2,1})$.
If $p=4$ then
the result is the four Fuchsian systems 
specified by the spectral  types
$\{(1,1)^5\}$,
$\{(2,1)^2,(1,1,1)^2 \}$,
$\{(3,1),(2,2)^2,(1,1,1,1)\}$,
and
$\{(2,2)^3,(2,1,1)\}$.
The first one has two deformation parameters and 
it corresponds 
to the Garnier system in two variables $(={\cal G}_{2,2})$.
The other three cases 
produce nonlinear ordinary differential equations 
of fourth order,
which have been investigated by
Sakai \cite{sak08} as
candidates of the master equations,
like $P_{\rm VI}$,
among the family of
 fourth-order Painlev\'e equations;
he clarified the polynomial Hamiltonian structure
and coalescence diagram
for each.
Note that the first and second 
of the three are equivalent
respectively to 
${\cal G}_{3,1}$ 
(see Example~\ref{ex:n=1}) 
and
to the fourth-order
Painlev\'e equation of type 
$D_6^{(1)}$
introduced by Sasano \cite{sas06}
(see also \cite{fs08}).
\end{remark}

\section{Polynomial Hamiltonian structure}
\label{sect:ham}

In this section
we present 
Hamiltonian formalism for the 
system ${\cal G}_{L,N}$
such that Hamiltonian functions are polynomials in the canonical variables.

The {\it Schlesinger system} 
is the following system of nonlinear differential equations
(see \cite{jmu81, sch12}):
\begin{equation}  \label{eq:ss}
\frac{\partial A_i}{\partial u_i}
= -\sum_{j \neq i} \frac{[A_i,A_j]}{u_i-u_j}, \quad
\frac{\partial A_i}{\partial u_j}
=\frac{[A_i,A_j]}{u_i-u_j} \quad (i \neq j)
\end{equation}
for $L \times L$ matrix-valued 
unknown functions $A_i$,
which describes the monodromy preserving deformations of 
a Fuchsian system of the form (\ref{eq:scf}).
Needless to say, 
${\cal G}_{L,N}$
is equivalent to a particular case of
the Schlesinger systems 
specified by the spectral type (\ref{eq:spec}).

Recall first that (\ref{eq:ss})
can be written as a Hamiltonian system 
(see, e.g., \cite{man99})
\[ \frac{\partial A_i}{\partial u_j}= \{ A_i , K_j\}
\]
with the Hamiltonian functions
\begin{equation}
K_i=  \frac{1}{2} 
\, \underset{z=u_i}{\rm{Res}}  \ {\rm tr} \, A^2
= \sum_{j \neq i} \frac{ {\rm tr} (A_iA_j) }{u_i-u_j},
\end{equation}
where the Poisson bracket $\{\mbox{\ },\mbox{\ }\}$
is given in a standard way by
\begin{equation} \label{eq:poisson_A}
\{ (A_i)_{m,n},(A_j)_{m',n'}\}
=\delta_{i,j}
\left( 
\delta_{m,n'}(A_i)_{m',n}  
-\delta_{m',n} (A_j)_{m,n'} 
\right).
\end{equation}
Moreover, 
a method to construct canonical variables for the above
Hamiltonian system  
has been established;
see \cite[Appendix 5]{jmms}.
Set $A_i=B^{(i)} C^{(i)}$ and define a Poisson bracket 
$\{\mbox{\ },\mbox{\ }\}$
over the space of matrices
$B^{(i)}$ and $C^{(i)}$ by
\[
\left\{
\left(B^{(i)}\right)_{m,n}, 
\left(C^{(i)}\right)_{n,m} 
\right\}=1  \quad
 \text{and}  \quad
\{\text{otherwise}\}=0.
\]
This Poisson bracket coincides with 
the previous one (\ref{eq:poisson_A}),
in fact.
Hence 
the Schlesinger system is equivalent to
the canonical Hamiltonian system 
attached
with the fundamental $2$-form
\[
\Gamma= \sum_{i=0}^{N+1} {\rm tr} \left( 
{\rm d} C^{(i)}  \wedge {\rm d} B^{(i)} \right)
- \sum_{i=0}^{N+1} {\rm d}  K_i \wedge  {\rm d} u_i .
\]
However, the above choice of canonical variables is 
redundant 
because it is possible to reduce
 the number of canonical variables 
 to that of accessory parameters of 
the Fuchsian system (\ref{eq:scf}).

Next we shall consider the Hamiltonian formalism of 
${\cal G}_{L,N}$
and carry out
the reduction of canonical variables.
In this case the fundamental $2$-form 
reads (see Sect.~\ref{sect:lax})
\begin{equation}  \label{eq:2-form}
\Gamma
=\sum_{i=0}^{N} \sum_{n=0}^{L-1} 
{\rm d} c_n^{(i)}  \wedge {\rm d} b_n^{(i)}
-\sum_{i=1}^{N} {\rm d}  K_i \wedge  {\rm d} u_i
\end{equation}
with 
\[
b_n^{(i)}
=\frac{ -g_{n,-n}^{(i)} }{  L {t_i}^n \prod_{m=1}^n f_{ m,-m+1 }^{(i)} }
\quad \text{and} \quad
c_n^{(i)}
={t_i}^n\prod_{m=1}^n f_{ m,-m+1 }^{(i)}
\]
for $0 \leq i \leq N$ and $0 \leq n \leq L-1$.
Here we have fixed 
$t_0=1$ and thereby $u_0=1$.
Observe that
\begin{equation} \label{eq:bc=}
b_n^{(0)} c_n^{(0)}
= - \kappa_{n} - \sum_{i=1}^N b_n^{(i)} c_n^{(i)},
\end{equation}
which follows from (\ref{eq:kappa_n}) by means of 
$b_n^{(i)} c_n^{(i)}= - g_{n,-n}^{(i)}/L$.
Accordingly the first term of (\ref{eq:2-form}) 
can be computed 
as follows:
\begin{align*}
\sum_{i=0}^{N} \sum_{n=0}^{L-1} 
{\rm d} c_n^{(i)}  \wedge {\rm d} b_n^{(i)} 
&=
\sum_{i=0}^{N} \sum_{n=1}^{L-1} 
{\rm d} c_n^{(i)}  \wedge {\rm d} b_n^{(i)},
 \quad \text{since $c_0^{(i)}=1$,}
\\
&=\sum_{i=0}^{N} \sum_{n=1}^{L-1} 
{\rm d} \log c_n^{(i)}  \wedge {\rm d} \left(b_n^{(i)} c_n^{(i)}\right)
\\
&
=\sum_{i=1}^{N} \sum_{n=1}^{L-1} 
{\rm d} \log \frac{c_n^{(i)}}{c_n^{(0)}}  
\wedge {\rm d} \left(b_n^{(i)} c_n^{(i)}\right), 
\quad \text{using (\ref{eq:bc=})},
\\
&=\sum_{i=1}^{N} \sum_{n=1}^{L-1} 
{\rm d} \left( \frac{c_n^{(i)}}{c_n^{(0)}}  \right)
\wedge {\rm d} \left(b_n^{(i)} c_n^{(0)}\right)
\\
&=\sum_{i=1}^{N} \sum_{n=1}^{L-1} 
 {\rm d} \left(-b_n^{(i)} c_n^{(0)}\right)
 \wedge  {\rm d} \left( \frac{c_n^{(i)}}{c_n^{(0)}}  \right).
\end{align*}
Let us now
introduce the canonical variables $q_n^{(i)}$ and $p_n^{(i)}$
$(1 \leq i \leq N; 1 \leq n \leq L-1)$
defined by
\begin{equation} \label{eq:can}
q_n^{(i)} =\frac{c_n^{(i)}}{c_n^{(0)}}, \quad
p_n^{(i)}=-b_n^{(i)} c_n^{(0)}, 
\end{equation}
whose number,
$2N(L-1)$,
 is just enough 
for the Hamiltonian
system under consideration;
see Remark~\ref{remark:scf}.
In addition we
take the change of independent variables 
\[
s_i=\frac{1}{u_i}={t_i}^{L}
\]
so that the resulting Hamiltonian function
\[
H_i=-\frac{K_i}{{s_i}^2}
=-\frac{{\rm tr}(A_iA_{N+1})}{s_i}+ \sum^N_{
\begin{subarray}{c} 
j=0  \\
j \neq i
\end{subarray}
}\frac{s_j {\rm  tr}(A_iA_j)}{s_i(s_i-s_j)} 
\]
becomes identical with the standard one of $P_{\rm VI}$ when $(L,N)=(2,1)$;
see Example \ref{ex:n=1}.
The fundamental $2$-form is then rewritten as
\[
\Gamma= \sum_{i=1}^N
\left(\sum_{n=1}^{L-1} {\rm d} p_n^{(i)} \wedge {\rm d} q_n^{(i)}
-{\rm d}  H_i \wedge  {\rm d}  s_i \right).
\]
For convenience
we extendedly use the symbols
$q_n^{(i)}$ and $p_n^{(i)}$
also for $i=0$ or $n=0$;
namely, we put
\begin{align} \label{eq:qp_add}
q_n^{(0)}
&=1,
\quad
p_n^{(0)}
\left(=-b_n^{(0)}c_n^{(0)} \right)
=
\kappa_{n} - \sum_{i=1}^{N} q_n^{(i)} p_n^{(i)},
\\
q_0^{(i)}
&=1, \quad 
p_0^{(i)}
\left(=-b_0^{(i)}c_0^{(i)} \right)
=\theta_i-\sum_{n=1}^{L-1} q_n^{(i)} p_n^{(i)},
\nonumber
\end{align}
by taking 
(\ref{eq:A=bc}) and (\ref{eq:bc=})
into account.
We have then the

\begin{lemma} \label{lemma:trace}
It holds that
\begin{align}
{\rm tr} (A_i A_j) 
&=  \sum_{m,n=0}^{L-1}
  q_m^{(i)}   p_m^{(j)} q_n^{(j)} p_n^{(i)},
  \label{eq:trAiAj}\\
{\rm tr}(A_i A_{N+1})  
&= -\sum_{n=0}^{L-1} e_n q_n^{(i)} p_n^{(i)} -\sum_{j=0}^N \sum_{0 \leq m <n \leq L-1}
  q_m^{(i)}   p_m^{(j)} q_n^{(j)} p_n^{(i)}
 \label{eq:trAiAN+1}
\end{align}  
for $i,j=0,1,\ldots,N$.
\end{lemma}

\pf
It follows from
$(A_i)_{m,n}= b_m^{(i)} c_n^{(i)}$
that ${\rm tr}(A_i A_j)
=\sum_{m,n=0}^{L-1}(A_i)_{n,m}(A_j)_{m,n}
=\sum_{m,n=0}^{L-1} b_n^{(i)} c_m^{(i)} b_m^{(j)} c_n^{(j)}$,
which thus yields (\ref{eq:trAiAj}) via
$b_n^{(i)}c_n^{(j)}=-p_n^{(i)}q_n^{(j)}$.
The diagonal entries of $A_iA_{N+1}$ read 
\[
\begin{array}{rr}
(0,0): & \left( e_0b_0^{(i)}+w_{0,1} b_1^{(i)}+w_{0,2} b_2^{(i)}+ \cdots + w_{0,L-1} b_{L-1}^{(i)} \right) c_0^{(i)},
\\
(1,1):& \left( e_1b_1^{(i)}+w_{1,2} b_2^{(i)}+ \cdots + w_{1,L-1} b_{L-1}^{(i)} \right) c_1^{(i)},
\\
\vdots \quad
&\vdots \qquad
\\
(L-2,L-2): & \left(e_{L-2}b_{L-2}^{(i)} +  w_{L-2,L-1}  b_{L-1}^{(i)}  \right)
c_{L-2}^{(i)},
\\
(L-1,L-1): & e_{L-1}b_{L-1}^{(i)}c_{L-1}^{(i)}.
\end{array}
\]
Therefore
${\rm tr}(A_iA_{N+1})=\sum_{n=0}^{L-1}e_n b_n^{(i)}c_n^{(i)}+\sum_{0 \leq m < n \leq L-1} w_{m,n} b_n^{(i)}  c_{m}^{(i)}$.
If we remember
$w_{m,n}=\sum_{j=0}^Nv_{m,n}^{(j)}=-\sum_{j=0}^N b_m^{(j)} c_n^{(j)}$,
then we find 
${\rm tr}(A_iA_{N+1})=\sum_{n=0}^{L-1}e_n b_n^{(i)}c_n^{(i)}-\sum_{j=0}^N\sum_{0 \leq m < n \leq L-1}  b_m^{(j)} c_n^{(j)}b_n^{(i)}  c_{m}^{(i)}$;
thus (\ref{eq:trAiAN+1})
is verified.
 \qed
\\

By virtue of Lemma~\ref{lemma:trace} together with (\ref{eq:qp_add}),
the Hamiltonian function $H_i$
can be explicitly expressed
as a polynomial in the $2N(L-1)$ canonical variables
$q_n^{(i)}$ and $p_n^{(i)}$
$(1 \leq i \leq N; 1 \leq n \leq L-1)$.
Finally we arrive at the

\begin{thm} \label{thm:caneq}
The system ${\cal G}_{L,N}$
is equivalent to the canonical Hamiltonian system
\begin{equation} \label{eq:caneq}
\frac{\partial q_n^{(i)}}{\partial s_j}=\frac{\partial H_j}{\partial p_n^{(i)}},
\quad 
\frac{\partial p_n^{(i)}}{\partial s_j}=-\frac{\partial H_j}{\partial q_n^{(i)}}
\quad 
\left(
\begin{array}{l}
  i,j =1, \ldots, N\\
 n=1, \ldots,   L-1
\end{array}
\right)
\end{equation}
where the Hamiltonian function $H_i$ is defined by
\begin{equation} \label{eq:caneq_ham}
s_i H_i=
\sum_{n=0}^{L-1} e_n q_n^{(i)} p_n^{(i)} 
+\sum_{j=0}^N \sum_{0 \leq m <n \leq L-1}
 q_m^{(i)}   p_m^{(j)} q_n^{(j)} p_n^{(i)} 
+\sum^N_{
\begin{subarray}{c} 
j=0  \\
j \neq i
\end{subarray}
}
\frac{s_j}{s_i-s_j} 
 \sum_{m,n=0}^{L-1}
 q_m^{(i)}   p_m^{(j)} q_n^{(j)} p_n^{(i)}
\end{equation}
and
\[
s_0=q_n^{(0)}
=q_0^{(i)}=1,
\quad
p_n^{(0)}
=
\kappa_{n} - \sum_{i=1}^{N} q_n^{(i)} p_n^{(i)},
\quad 
p_0^{(i)}
=\theta_i-\sum_{n=1}^{L-1} q_n^{(i)} p_n^{(i)}.
\]
\end{thm}

We write the constant parameters 
contained in (\ref{eq:caneq})
as 
\begin{equation} \label{eq:vec_kappa}
\vec{\kappa}=(e_0,e_1,\ldots,e_{L-1},\kappa_0,\kappa_1,\ldots,\kappa_{L-1},\theta_0,\theta_1,\ldots,\theta_N),
\end{equation}
whose number is essentially $2 L +N-1$ 
according to (\ref{eq:rel_exp}).
Let ${\cal H}_{L,N}={\cal H}_{L,N}(\vec{\kappa})$
denote the Hamiltonian system (\ref{eq:caneq}).
Since all the differential equations 
originate from 
a single equation (\ref{eq:UCH}),
the system ${\cal H}_{L,N}$
is {\it a priori} completely integrable 
(in the Frobenius sense).
Or it can be shown
directly
by noticing the following facts:
(i) the $1$-form
$\omega = \sum_{i=1}^N H_i {\rm d}s_i$ is closed
for an arbitrary solution of (\ref{eq:caneq});
(ii) the relation
\[
\left(\frac{ \partial}{\partial s_j}\right) H_i = \left(\frac{ \partial}{\partial s_i}\right) H_j
= \frac{1}{(s_i-s_j)^2} 
 \sum_{m,n=0}^{L-1}
 q_m^{(i)}   p_m^{(j)} q_n^{(j)} p_n^{(i)}
 \quad (i \neq j)
\]
holds,
where the symbol $(\partial/\partial s_i)$
denotes the differentiation 
such that $q_n^{(i)}$ and $p_n^{(i)}$ are viewed to be 
independent of $s_i$.
These facts imply the commutativity of the flows induced by 
$H_1,H_2, \ldots, H_N$.

The correspondence between the
canonical variables $q_n^{(i)}$ and $p_n^{(i)}$
and the dependent variables  
given in Sect.~\ref{sect:sim}
is summarized as 
\begin{subequations} \label{subeq:qp_tau}
\begin{align} 
q_n^{(i)}
&= 
\frac{c_n^{(i)}}{c_n^{(0)}}
=
\left(\frac{t_i}{t_0}\right)^n 
\prod_{m=1}^n \frac{f_{m,-m+1}^{(i)}}{f_{m,-m+1}^{(0)} }
= \left(\frac{t_i}{t_0}\right)^n
\frac{ \sigma_{n,-n}(\theta_i+1) \sigma_{0,0}(\theta_0+1)  }{ \sigma_{0,0}(\theta_i+1) \sigma_{n,-n}(\theta_0+1) },
\label{eq:qp_tau_1}
\\
 q_n^{(i)} p_n^{(i)} 
 &=-b_n^{(i)}c_n^{(i)}
 = \frac{ g_{n,-n}^{(i)}}{L}= 
\frac{\theta_i}{L}
\frac{\sigma_{n-1,-n-1} (\theta_i-1) \sigma_{n,-n}(\theta_i+1)
}{   \sigma_{n,-n-1} \sigma_{n-1,-n} }.
\label{eq:qp_tau_2}
\end{align}
\end{subequations}

\begin{example}[Case $N=1$]\rm
\label{ex:n=1}
Let us restrict ourselves to
 the case $N=1$;
thus
${\cal H}_{L,1}$ 
becomes a system
of ordinary differential equations.
We begin with the case $L=2$, 
which is the first nontrivial one.
Write 
$(q,p,H,s)= (q_1^{(1)},p_1^{(1)},H_1,s_1)$
and $\theta=\theta_1$.
Then the Hamiltonian function 
can be expressed as 
\[H=H_{\rm VI}(a_0,a_1,a_2,a_3,a_4;q,p)
+ 
\frac{ \theta (e_0 (s-1)+\kappa_0-\theta) }{ s(s-1)}
\]
under the substitution 
\[
a_0= e_0-e_1+ \kappa_1+1, \quad
a_1= -\kappa_1+\theta,
\quad
a_2=-\theta,
\quad a_3= -e_0+e_1+\kappa_0,
\quad
a_4=-\kappa_{0}+\theta.
\]
Here 
$H_{\rm VI}=H_{\rm VI}(a_0,a_1,a_2,a_3,a_4;q,p)$
denotes the Hamiltonian function 
of $P_{\rm VI}$
and is defined by
\begin{align*}
s(s-1) H_{\rm VI}
&= q(q-1)(q-s)p^2
\\
& \quad
-\left(
(a_0-1)q(q-1) +
a_3 q(q-s)
+ a_4 (q-1)(q-s)
\right) p
\\
&\quad
+ a_2(a_1+a_2) q
\end{align*}
with $a_i$ being constant parameters such that
$a_0+a_1+2a_2+a_3+a_4=1$;
see \cite{malm22, oka87}.

Now we turn to the case of general $L \geq 2$.
Let $(q_n,p_n,H,s)= (q_n^{(1)},p_n^{(1)},H_1,s_1 )$
and $\theta=\theta_1$.
Then the Hamiltonian function
of ${\cal H}_{L,1}$
 takes a
{\it coupled} form of $P_{\rm VI}$ ones 
as follows:
\begin{align} \label{eq:cp6}
H&= 
\sum_{n=1}^{L-1} 
H_{\rm VI}(a_{0,n},a_{1,n},a_{2,n},a_{3,n},a_{4,n};q_n,p_n)
+ 
\frac{ \theta (e_0 (s-1)+\kappa_0-\theta) }{ s(s-1)}
 \\
& \quad
+ \sum_{1 \leq m < n \leq L-1}
\frac{(q_m-1)p_mq_n((q_n-s)p_n-\kappa_n)
+(q_n-s)p_nq_m((q_m-1)p_m-\kappa_m)}{s(s-1)},
\nonumber
\end{align}
where  
the last term reflects an interaction and
the correspondence of constant parameters 
reads
\[
a_{0,n}= e_0-e_n+ \kappa_n+1, \quad
a_{1,n}= -\kappa_n+\theta,
\quad
a_{2,n}=-\theta,
\quad 
a_{3,n}= -e_0+e_n+\kappa_0,
\quad
a_{4,n}=-\kappa_{0}+\theta.
\]
Interestingly enough,
as has been pointed out
by Fuji and Suzuki (see \cite{fs09, suz10}),
the coupled Hamiltonian 
(\ref{eq:cp6})
can be derived alternatively from the deformation of a certain linear system 
that is not Fuchsian but has one regular and one irregular singularities;
cf. Sect.~\ref{sect:lax}.
It is expected to exist some integral transform
(like a Laplace one)
between the two kinds of Lax formalism. 
\end{example}

\begin{remark} \rm
We cite the recent result by Dubrovin and Mazzocco
\cite{dm07};
they have studied Hamiltonian formalism
of the Schlesinger system associated with 
the general spectral type
(cf. (\ref{eq:spec})).
Their construction is based on 
a scalar differential equation of higher order
that is reduced from a Fuchsian system of the form
(\ref{eq:scf});
and the apparent singularities 
(see \cite{ko83})
produced by the reduction procedure
are adopted as
the half of the canonical variables, 
i.e., the generalized coordinates.
%
The resulting Hamiltonian functions are rational in the canonical variables.
It would be an interesting problem to transform the 
{\it general} Schlesinger system into a Hamiltonian system
with the Painlev\'e property (see \cite{malg83, miw81})
whose
 Hamiltonian functions are polynomials in the canonical variables,
 like ${\cal H}_{L,N}$.
\end{remark}

\allowdisplaybreaks

\section{Birational canonical transformations}
\label{sect:weyl}
This section
is devoted to birational symmetries of the Hamiltonian system
${\cal H}_{L,N}={\cal H}_{L,N}(\vec{\kappa})$.
Here, to be precise, a birational canonical transformation
of variables $(q_n^{(i)},p_n^{(i)},s_i)$ is said to be a {\it symmetry} if
it keeps the system invariant except changing the constant
parameters
$\vec{\kappa}$.

First we translate the action of
$\langle r_k \rangle$ 
discussed in Sect.~\ref{sect:hom}
into birational canonical transformations of ${\cal H}_{L,N}$.
Note that
$\langle r_k \rangle$ is isomorphic to 
an affine Weyl group of type 
$A_{L-1}^{(1)}$, 
denoted by $W(A_{L-1}^{(1)})$.
For each $k \in {\mathbb Z}/L{\mathbb Z}$, 
let $r_k(\sigma_{k,n})=\hat{\sigma}_{k,n}$ and  
$r_k(\sigma_{m,n})=\sigma_{m,n}$ $(m \neq k)$.
Substitution of (\ref{eq:special_gar})
and $t=t_i$
in (\ref{eq:diff_hat_1})
yields
\[
D_i
\hat{\sigma}_{k,n} \cdot \sigma_{k,n} 
= \theta_i \sigma_{k-1,n}(\theta_i-1) \sigma_{k+1,n}(\theta_i+1)
\]
with the aid of (\ref{eq:tdt}).
Therefore we have
\[
\sum_{i=0}^N
t_i D_i
\hat{\sigma}_{k,n} \cdot \sigma_{k,n} 
= \sum_{i=0}^N\theta_i t_i \sigma_{k-1,n}(\theta_i-1) \sigma_{k+1,n}(\theta_i+1).
\]
In view of the homogeneity (\ref{eq:hom_gar})
we conclude that 
\begin{equation} \label{eq:weyl_A_1}
\hat{\sigma}_{k,n} = \frac{1}{\alpha_k \sigma_{k,n}} \sum_{i=0}^N\theta_i t_i \sigma_{k-1,n}(\theta_i-1) \sigma_{k+1,n}(\theta_i+1); 
\end{equation}
recall  (\ref{eq:alpha}).
Similarly, we deduce from (\ref{eq:diff_hat_2}) that
\begin{equation}  \label{eq:weyl_A_2}
t_i \sigma_{k-1,n+1}(\theta_i-1) 
\sigma_{k+1,n}
- \hat{\sigma}_{k,n+1}
 \sigma_{k,n}(\theta_i-1) 
  +\hat{\sigma}_{k,n}(\theta_i-1) 
  \sigma_{k,n+1}=0.
\end{equation}
Through (\ref{eq:gar_f}) and (\ref{eq:gar_g}),
the action of $r_k$ on $(f_{m,n}^{(i)},g_{m,n}^{(i)})$ 
is determined by
(\ref{eq:weyl_A_1}) and (\ref{eq:weyl_A_2})
as follows:
\begin{align*}
r_k(f_{k,n}^{(i)})
& =
f_{k,n}^{(i)} \left(
1+ \frac{ \alpha_k t_i f_{k+1,n-1}^{(i)}}{\sum_{j=0}^N  t_j f_{k+1,n-1}^{(j)} g_{k,n-1}^{(j)}}
\right),
\\
r_k(f_{k+1,n-1}^{(i)})
&=f_{k+1,n-1}^{(i)}  \left(
1- \frac { \alpha_k t_i f_{k+1,n-1}^{(i)}}{  \alpha_k t_i f_{k+1,n-1}^{(i)} + \sum_{j=0}^N  t_j f_{k+1,n-1}^{(j)} g_{k,n-1}^{(j)} }
\right),
\\
r_k(g_{k,n}^{(i)})
&=g_{k,n}^{(i)} \left(
 1+ \frac{ \alpha_k t_i f_{k+1,n}^{(i)}}{\sum_{j=0}^N  t_j f_{k+1,n}^{(j)} g_{k,n}^{(j)}}
 \right),
 \\
 r_k(g_{k+1,n-1}^{(i)})
 &=g_{k+1,n-1}^{(i)}-
\frac{ \alpha_k t_i f_{k+1,n}^{(i)} g_{k,n}^{(i)} }{\sum_{j=0}^N  t_j f_{k+1,n}^{(j)} g_{k,n}^{(j)}},
\end{align*}
for $n \in {\mathbb Z}/L{\mathbb Z}$.
It is then easy to construct  the corresponding
transformation of
$(q_n^{(i)},p_n^{(i)})$
by virtue of (\ref{subeq:qp_tau}).
Moreover,
as has been mentioned in Sect.~\ref{sect:hom},
the system ${\cal H}_{L,N}$
enjoys another action $\langle {r_k}'\rangle$
of $W(A_{L-1}^{(1)})$
associated with the root variables
$\beta_k=-d_{m,k-1}+ 2 d_{m,k}-d_{m,k+1} +1$,
which commutes with
 the previous one $\langle {r_k}\rangle$.

Next
we observe that a cyclic permutation of the suffixes 
$\pi: (\sigma_{m,n},d_{m,n}) \mapsto  (\sigma_{m+1,n-1},d_{m+1,n-1})$  
keeps the bilinear expression
(\ref{subeq:bil_gar})
of ${\cal H}_{L,N}$
invariant,
and so does the interchange of suffixes 
$\rho: (\sigma_{m,n}, d_{m,n}, t_i) \mapsto  (\sigma_{n,m}, -d_{n,m}, 1/t_i)$.
These trivial symmetries can be lifted to birational canonical transformations of ${\cal H}_{L,N}$.
Note that $\pi$ realizes a Dynkin automorphism 
which rotates simultaneously the two Dynkin diagrams of type $A_{L-1}^{(1)}$
and that
$\rho$ represents an interchange of the two diagrams.

For notational simplicity we extend the suffix $n$ of the canonical variables 
$(q_n^{(i)},p_n^{(i)})$
and parameters
$e_n$  and $\kappa_n$ for any $n \in {\mathbb Z}$
by the conditions
(cf. (\ref{subeq:qp_tau}))
\[
q_{n+L}^{(i)}=s_i q_{n}^{(i)}, \quad 
p_{n+L}^{(i)}= \frac{p_{n}^{(i)}}{s_i}, \quad
e_{n+L}=e_n+1, \quad \kappa_{n+L}=\kappa_n.
\]
We set
\begin{equation} \label{eq:ab}
\mathfrak{a}_n=\frac{\alpha_n}{L}=e_{n+1}-e_n,
\quad 
\mathfrak{b}_n=\frac{\beta_n}{L}=e_{L-n}-e_{L-n-1}-\kappa_{L-n}+\kappa_{L-n-1}
\end{equation}
for $0 \leq n \leq L-1$.
It thus holds that
$\sum_{n=0}^{L-1}\mathfrak{a}_n=\sum_{n=0}^{L-1}\mathfrak{b}_n=1$.

We now state the result.

\begin{thm}  \label{thm:symmetry_1} 
The Hamiltonian system
${\cal H}_{L,N}(\vec{\kappa})$
is invariant under the birational canonical transformations
$r_n$, ${r_n}'$, $\pi$, and $\rho$
$(n = 0,1,\ldots,L-1 )$
defined as follows{\rm:}
\begin{itemize}
\item Action on the parameters $\vec{\kappa}$.
\begin{align*}
r_n&:e_{n} \mapsto e_n + {\mathfrak a}_n, \quad 
e_{n+1} \mapsto e_{n+1} - {\mathfrak a}_n, \quad
\kappa_{n} \mapsto \kappa_n + {\mathfrak a}_n, \quad 
\kappa_{n+1} \mapsto \kappa_{n+1} - {\mathfrak a}_n.
\\
{r_n}'&:\kappa_{L-n} \mapsto \kappa_{L-n} + {\mathfrak b}_n, \quad 
\kappa_{L-n-1} \mapsto \kappa_{L-n-1} - {\mathfrak b}_n.
\\
\pi&:e_n \mapsto e_{n+1}-\frac{1}{L}, \quad 
\kappa_n \mapsto \kappa_{n+1}.
\\
\rho &:e_n \mapsto \kappa_{L-n}-e_{L-n}-\frac{\sum_{i=0}^N \theta_i }{L} +1,
\quad
\kappa_n \mapsto \kappa_{L-n}.
\end{align*}
\item Action on the canonical variables $(q_n^{(i)},p_n^{(i)})$.
\begin{align*}
r_n \ (n \neq 0)
&: 
\left\{
\begin{array}{l} \displaystyle
q_n^{(i)} \mapsto q_n^{(i)} + \frac{{\mathfrak a}_n (q_{n+1}^{(i)}-q_n^{(i)})}{ {\mathfrak a}_n 
+ \sum_{j=0}^N q_{n+1}^{(j)} p_n^{(j)}},
\\
\displaystyle
p_n^{(i)} \mapsto 
p_n^{(i)} \left(
1+ \frac{{\mathfrak a}_n}{\sum_{j=0}^N q_{n+1}^{(j)} p_n^{(j)}} \right), 
\\
\displaystyle
p_{n+1}^{(i)}
\mapsto
p_{n+1}^{(i)}-  
\frac{{\mathfrak a}_n p_n^{(i)}}{\sum_{j=0}^N q_{n+1}^{(j)} p_n^{(j)}}.
\end{array}
\right.
\\
r_0 
&:
\left\{
\begin{array}{l} \displaystyle
q_n^{(i)} \mapsto 
q_n^{(i)} \left(  1- \frac{ {\mathfrak a}_0 (q_1^{(i)} - 1) }{ {\mathfrak a}_0 q_1^{(i)} + \sum_{j=0}^N q_1^{(j)} p_0^{(j)}} \right),
\\
\displaystyle
p_n^{(i)} \mapsto
p_n^{(i)} \left(  1+ \frac{ {\mathfrak a}_0 (q_1^{(i)} - 1) }{ {\mathfrak a}_0  + \sum_{j=0}^N q_1^{(j)} p_0^{(j)}} \right)
\quad (n \neq 1),
\\
\displaystyle
p_1^{(i)} \mapsto
\left(
p_1^{(i)}-\frac{ {\mathfrak a}_0 p_0^{(i)}}{ \sum_{j=0}^N q_1^{(j)} p_0^{(j)} }
\right)
 \left(  1+ \frac{ {\mathfrak a}_0 (q_1^{(i)} - 1) }{ {\mathfrak a}_0  + \sum_{j=0}^N q_1^{(j)} p_0^{(j)}} \right).
\end{array}
\right.
\\
{r_n}' \ (n \neq 0)
&: 
\left\{
\begin{array}{l} \displaystyle
q_{L-n}^{(i)} \mapsto q_{L-n}^{(i)} + \frac{{\mathfrak b}_n (q_{L-n-1}^{(i)}-q_{L-n}^{(i)})}{ {\mathfrak b}_n 
+ \sum_{j=0}^N q_{L-n-1}^{(j)} p_{L-n}^{(j)}},
\\
\displaystyle
p_{L-n}^{(i)} \mapsto 
p_{L-n}^{(i)} \left(
1+ \frac{{\mathfrak b}_n}{\sum_{j=0}^N q_{L-n-1}^{(j)} p_{L-n}^{(j)}} \right),
\\
\displaystyle
p_{L-n-1}^{(i)}
\mapsto
p_{L-n-1}^{(i)}-  
\frac{{\mathfrak b}_n p_{L-n}^{(i)}}{\sum_{j=0}^N q_{L-n-1}^{(j)} p_{L-n}^{(j)}}.
\end{array}
\right.
\\
{r_0}'
&:
\left\{
\begin{array}{l} \displaystyle
q_n^{(i)} \mapsto 
q_n^{(i)} \left(  1- \frac{ {\mathfrak b}_0 (q_{-1}^{(i)} - 1) }{ {\mathfrak b}_0 q_{-1}^{(i)} + \sum_{j=0}^N q_{-1}^{(j)} p_0^{(j)}} \right),
\\
\displaystyle
p_n^{(i)} \mapsto
p_n^{(i)} \left(  1+ \frac{ {\mathfrak b}_0 (q_{-1}^{(i)} - 1) }{ {\mathfrak b}_0   + \sum_{j=0}^N q_{-1}^{(j)} p_0^{(j)}} \right)
\quad (n \neq L-1), 
\\
\displaystyle
p_{L-1}^{(i)} \mapsto
\frac{1}{s_i}
\left(
p_{-1}^{(i)}-\frac{ {\mathfrak b}_0 p_0^{(i)}}{ \sum_{j=0}^N q_{-1}^{(j)} p_0^{(j)} }
\right)
 \left(  1+ \frac{ {\mathfrak b}_0 (q_{-1}^{(i)} - 1) }{ {\mathfrak b}_0  + \sum_{j=0}^N q_{-1}^{(j)} p_0^{(j)}} \right).
\end{array}
\right.
\\
\pi
&:
q_n^{(i)} \mapsto 
\frac{q_{n+1}^{(i)}}{q_1^{(i)}}, \quad 
p_n^{(i)} \mapsto  p_{n+1}^{(i)} q_1^{(i)}.
\\
\rho 
&:
s_i \mapsto  \frac{1}{s_i}, \quad 
q_n^{(i)} \mapsto
\frac{q_{L-n}^{(i)}}{s_i}, \quad 
p_n^{(i)} \mapsto s_i p_{L-n}^{(i)}. 
\end{align*}
\end{itemize}
{\rm(}Here we have omitted to write the action on the variables if it is trivial.{\rm)}
Moreover, 
these transformations satisfy the relations{\rm:}
${r_n}^2=(r_n r_{n \pm 1})^3=({r_n}')^2=({r_n}' {r_{n \pm 1}}')^3=\pi^L=\rho^2=\id$,
$\pi r_{n} =r_{n+1} \pi $, 
$\pi {r_n}'  = {r_{n-1}}' \pi$,
and 
$ \rho  r_n={r_n}' \rho$.
\end{thm}

Let us explore further symmetries of 
${\cal H}_{L,N}$
besides those in Theorem~\ref{thm:symmetry_1}.
First we consider a symmetry shifting the parameter $\theta_i$ to $\theta_i -1$
at the level of the variables $f_{m,n}^{(i)}$ and $g_{m,n}^{(i)}$.
It readily follows from (\ref{eq:gar_f}) and  (\ref{eq:gar_g})
that
\begin{equation} \label{eq:Ti_1}
f_{m,n}^{(i)}(\theta_i-1)= \frac{g_{m,n}^{(i)}}{ g_{m+1,n-1}^{(i)}}
f_{m+1,n}^{(i)}.
\end{equation}
Combining this with (\ref{eq:uv_fg_4}) shows that
 \begin{equation} \label{eq:Ti_2}
 f_{m,n}^{(j)}(\theta_i-1)= \frac{t_i U_{m,n}^{(i,j)}g_{m,n}^{(i)}}{t_j V_{m,n}^{(i,j)}g_{m+1,n-1}^{(i)}} f_{m+1,n}^{(i)}
 \quad (i \neq j).
 \end{equation}
We observe
for $i \neq j$
that
\begin{align}  \label{eq:Ti_3}
g_{m,n}^{(j)}(\theta_i-1) &= \theta_j \frac{\sigma_{m-1,n-1}(\theta_i-1,\theta_j-1)\sigma_{m,n}(\theta_i-1,\theta_j+1)}{\sigma_{m,n-1}(\theta_i-1)\sigma_{m-1,n}(\theta_i-1)}
 \\
&=\frac{t_i-t_j}{t_j}
\frac{\sigma_{m-1,n-1}(\theta_i-1,\theta_j-1)\sigma_{m,n}}{\sigma_{m,n-1}(\theta_i-1)\sigma_{m-1,n}(\theta_j-1)} 
\times
\frac{g_{m,n+1}^{(j)}}{g_{m,n+1}^{(i)}} U_{m,n+1}^{(i,j)}
\nonumber \\
&= \left( \frac{t_i}{t_j}
\frac{\sigma_{m,n-1}(\theta_j-1)\sigma_{m-1,n}(\theta_i-1) }{\sigma_{m,n-1}(\theta_i-1)\sigma_{m-1,n}(\theta_j-1)}
-1
\right) 
\frac{g_{m,n+1}^{(j)}}{g_{m,n+1}^{(i)}} U_{m,n+1}^{(i,j)},
\quad \text{using (\ref{eq:bil_gar_1})},
\nonumber \\
&= 
 \left( \frac{t_i f_{m,n+1}^{(i)}(\theta_i-1)}{ t_j f_{m,n+1}^{(j)}(\theta_j-1) }
 -1 \right)\frac{g_{m,n+1}^{(j)}}{g_{m,n+1}^{(i)}} U_{m,n+1}^{(i,j)},
\quad \text{using  (\ref{eq:gar_f})},
\nonumber \\
&= 
 \left( \frac{t_i f_{m+1,n+1}^{(i)} g_{m+1,n}^{(j)}}{ t_j f_{m+1,n+1}^{(j)} g_{m+1,n}^{(i)} }
 -\frac{g_{m,n+1}^{(j)}}{g_{m,n+1}^{(i)}} \right)U_{m,n+1}^{(i,j)},
\quad \text{using  (\ref{eq:Ti_1})}.
\nonumber
\end{align}  
By (\ref{eq:kappa}) we have
\begin{equation}  \label{eq:Ti_4}
g_{m,n}^{(i)}(\theta_i-1)= \kappa_{m,n}-1-\sum_{j \neq i} g_{m,n}^{(j)}(\theta_i-1).
\end{equation}
The transformations (\ref{eq:Ti_1})--(\ref{eq:Ti_4}) 
provide a symmetry of the system ${\cal G}_{L,N}$, 
(\ref{subeq:gar}), 
shifting the parameter $\theta_i$ to $\theta_i-1$;
however, they do not naively give a symmetry of ${\cal H}_{L,N}$.
To reach a birational canonical transformation of ${\cal H}_{L,N}$,
we need to combine  
a trivial symmetry of 
(\ref{subeq:gar})
shifting the suffixes:
$(f_{m,n}^{(i)},g_{m,n}^{(i)},d_{m,n}) \mapsto (f_{m-1,n}^{(i)},g_{m-1,n}^{(i)},d_{m-1,n})$.
As a result we obtain a symmetry $\eta_i$ of 
${\cal H}_{L,N}$
which acts on the parameters 
as $\theta_i \mapsto \theta_i-1$ and ${\mathfrak a}_n \mapsto {\mathfrak a}_{n-1}$; 
see Theorem~\ref{thm:symmetry_2} below.
We do not go into detail of computations.

It is easy to find a group of symmetries $\langle \zeta_{ij} \rangle (\simeq {\mathfrak S}_{N+1})$,
which 
is generated by a permutation 
of the singularities
$z=u_i=1/s_i$ 
$(0 \leq i \leq N)$
of 
the associated Fuchsian system;
see Sect.~\ref{sect:lax}.

Finally
we 
deal with a symmetry deduced from the bilinear expression of
${\cal H}_{L,N}$ again.
Observe that (\ref{subeq:bil_gar})
is invariant under the transformation
\[\iota:
\sigma_{m,n}=
\sigma_{m,n}({\boldsymbol \theta},{\boldsymbol t})
\mapsto 
\sigma_{-m-1,-n-1}(-{\boldsymbol \theta},{\boldsymbol t}),
\quad
d_{m,n} \mapsto  d_{-m-1,-n-1}, \quad
\theta_i \mapsto -\theta_i. \]
Hence we have
\begin{align}   \label{eq:iota_1}
\iota(q_n^{(i)}) &= \iota \left(  \left(\frac{t_i}{t_0}\right)^n
\frac{ \sigma_{n,-n}(\theta_i+1) \sigma_{0,0}(\theta_0+1)  }{ \sigma_{0,0}(\theta_i+1) \sigma_{n,-n}(\theta_0+1) }\right),
\quad \text{using (\ref{eq:qp_tau_1})},
\\
&=  \left(\frac{t_i}{t_0}\right)^n
\frac{ \sigma_{-n-1,n-1}(\theta_i-1) \sigma_{-1,-1}(\theta_0-1)  }{ \sigma_{-1,-1}(\theta_i-1) \sigma_{-n-1,n-1}(\theta_0-1) }
\nonumber \\
&= \left(\frac{t_i}{t_0}\right)^n \prod_{m=-n}^{-1} \frac{f_{m,-m}^{(i)}(\theta_i-1)}{f_{m,-m}^{(0)}(\theta_0-1)}
\nonumber \\
&=
 \left(\frac{t_i}{t_0}\right)^n \prod_{m=-n}^{-1}
  \frac{  g_{m,-m}^{(i)} g_{m+1,-m-1}^{(0)}  }{ g_{m,-m}^{(0)} g_{m+1,-m-1}^{(i)}}
  \frac{f_{m+1,-m}^{(i)}}{f_{m+1,-m}^{(0)}},
\quad
\text{using (\ref{eq:Ti_1})},
\nonumber \\
&= 
\frac{  g_{-n,n}^{(i)} g_{0,0}^{(0)}  }{ g_{-n,n}^{(0)} g_{0,0}^{(i)}}
 \frac{1}{q_{-n}^{(i)}} 
\nonumber \\
 &= \frac{s_ip_{L-n}^{(i)} p_0^{(0)}}{p_{L-n}^{(0)} p_{0}^{(i)}}.
\nonumber 
\end{align}
Similarly, it follows that 
\begin{align}  \label{eq:iota_2}
\iota(q_{n}^{(i)}p_n^{(i)})
&= \iota 
\left(\frac{\theta_i}{L}
\frac{\sigma_{n-1,-n-1} (\theta_i-1) \sigma_{n,-n}(\theta_i+1)
}{   \sigma_{n,-n-1} \sigma_{n-1,-n} }
\right),
\quad \text{using (\ref{eq:qp_tau_2})},
 \\
&=\frac{-\theta_i}{L}
\frac{\sigma_{-n,n} (\theta_i+1) \sigma_{-n-1,n-1}(\theta_i-1)
}{   \sigma_{-n-1,n} \sigma_{-n,n-1} }
\nonumber
\\
&= -q_{L-n}^{(i)} p_{L-n}^{(i)}.
\nonumber
\end{align}
These formulae (\ref{eq:iota_1}) and (\ref{eq:iota_2}) define a birational canonical transformation of ${\cal H}_{L,N}$.

The above results are summed up in the 
  
\begin{thm}  \label{thm:symmetry_2} 
The Hamiltonian system
${\cal H}_{L,N}(\vec{\kappa})$
is invariant under 
the birational canonical transformations
$\eta_i$, $\zeta_{ij}$, and $\iota$
$(i,j=0,1, \ldots, N; i \neq j)$
defined as follows{\rm:}
\begin{itemize}
\item Action on the parameters $\vec{\kappa}$.
\begin{align*}
\eta_i &: e_n \mapsto e_{n-1}+\frac{1}{L}, \quad
\kappa_n \mapsto \kappa_n-e_n+e_{n-1}, \quad
\theta_i \mapsto \theta_i-1.
\\
\zeta_{ij} &: \theta_i \leftrightarrow \theta_j. 
\\
\iota &: e_{n} \mapsto -e_{L-n}+1, \quad
\kappa_n \mapsto -\kappa_{L-n}, \quad 
\theta_i \mapsto - \theta_i.
\end{align*}
\item Action on the canonical variables $(q_n^{(i)},p_n^{(i)})$.
\end{itemize}
\begin{align*}
\eta_i&:
\left\{
\begin{array}{l}
\displaystyle
q_n^{(j)}
\mapsto
\eta_i(q_n^{(j)})
= \frac{ \left(\sum_{m=1}^{L} p_{-m}^{(i)}\right) \left(\sum_{m=1}^L p_{n-m}^{(i)} q_{n-m}^{(j)} \right)}{ \left(\sum_{m=1}^{L} p_{n-m}^{(i)}\right) \left(\sum_{m=1}^L p_{-m}^{(i)} q_{-m}^{(j)} \right) }
\quad (\text{for $\forall j$}),
\\
\displaystyle
p_n^{(j)}
\mapsto \eta_i(p_n^{(j)})=\frac{1}{\eta_i(q_n^{(j)})} \frac{s_j}{s_i-s_j} 
\left( \frac{p_n^{(j)}}{p_n^{(i)}}- \frac{p_{n-1}^{(j)}}{p_{n-1}^{(i)}}
\right) \sum_{m=1}^L p_{n-m}^{(i)} q_{n-m}^{(j)} 
\quad (j \neq i),
\\
\displaystyle
p_n^{(i)}
\mapsto 
\frac{1}{\eta_i(q_n^{(i)})}\left(
\kappa_n-e_n+e_{n-1}-\sum_{j \neq i} \eta_i(q_n^{(j)}p_n^{(j)})
\right).
\end{array}
\right. 
\\
\zeta_{ij} \ 
(i,j \neq 0)
&: s_i \leftrightarrow s_j, \quad 
 q_n^{(i)} \leftrightarrow q_n^{(j)}, \quad 
p_n^{(i)} \leftrightarrow p_n^{(j)}.
\\
\zeta_{i0}=\zeta_{0i}
&: 
\left\{
\begin{array}{l}
\displaystyle
s_i \mapsto \frac{1}{s_i}, \quad s_j \mapsto \frac{s_j}{s_i}, 
\quad
q_n^{(i)} \mapsto \frac{1}{q_n^{(i)}}, \quad 
q_n^{(j)} \mapsto \frac{q_n^{(j)}}{q_n^{(i)}},
\\
\displaystyle
p_n^{(i)} \mapsto q_n^{(i)}p_n^{(0)},
\quad p_n^{(j)} \mapsto q_n^{(i)}p_n^{(j)}
\quad (j \neq i).
\end{array}
\right.
\\
\iota &:
q_n^{(i)}
\mapsto
\frac{s_i p_{L-n}^{(i)} p_0^{(0)}}{p_{L-n}^{(0)} p_{0}^{(i)}},
\quad
p_n^{(i)} \mapsto 
- \frac{q_{L-n}^{(i)} p _0^{(i)} p_{L-n}^{(0)} }{s_i p_0^{(0)} }.
\end{align*}
\end{thm}

\begin{remark} \rm
We may regard
$\mathfrak{a}_n$, $\mathfrak{b}_n$ ($1 \leq n \leq L-1$)
and $\theta_i$ ($0 \leq i \leq  N$)
as the $2L+N-1$ constant parameters of ${\cal H}_{L,N}$
instead of $\vec{\kappa}$;
see (\ref{eq:vec_kappa}) and (\ref{eq:ab}).
For reference
we summarize
how
the birational symmetries in 
Theorems~\ref{thm:symmetry_1} and 
\ref{thm:symmetry_2} 
act on $\mathfrak{a}_n$, $\mathfrak{b}_n$,
and $\theta_i$ below.
\begin{align*}
r_n&:{\mathfrak a}_{n} \mapsto -{\mathfrak a}_n, \quad 
{\mathfrak a}_{n \pm 1} \mapsto {\mathfrak a}_{n \pm 1} + {\mathfrak a}_n.
\\
{r_n}'&:{\mathfrak b}_{n} \mapsto -{\mathfrak b}_n, \quad 
{\mathfrak b}_{n \pm 1} \mapsto {\mathfrak b}_{n \pm 1} + {\mathfrak b}_n.
\\
\pi&:{\mathfrak a}_{n} \mapsto {\mathfrak a}_{n+1}, \quad 
{\mathfrak b}_{n} \mapsto {\mathfrak b}_{n-1}.
\\
\rho &:{\mathfrak a}_{n} \leftrightarrow {\mathfrak b}_{n}.
\\
\eta_i &: 
{\mathfrak a}_{n} \mapsto {\mathfrak a}_{n-1}, \quad 
\theta_i \mapsto \theta_i-1.
\\
\zeta_{ij} &: \theta_i \leftrightarrow \theta_j. 
\\
\iota &: 
{\mathfrak a}_{n} \mapsto {\mathfrak a}_{L-n}, \quad 
{\mathfrak b}_{n} \mapsto {\mathfrak b}_{L-n}, \quad 
\theta_i \mapsto - \theta_i.
\end{align*}
Recall that
the groups of canonical transformations
$\langle r_n \rangle$ and $\langle {r_n}' \rangle$ mutually commute 
and each of them
gives a birational realization of $W(A_{L-1}^{(1)})$.
\end{remark}

\begin{remark}[Additional symmetry valid for only $N=1$]
\label{remark:add_sym}
\rm
Let us consider the case $N=1$.
Write $(q_n,p_n,s)=(q_n^{(1)},p_n^{(1)},s_1)$ and $\theta=\theta_1$.
We then find another symmetry
given as follows
(see also the appendix):
\[
\varphi: 
\left\{
 \begin{array}{l}
e_0 \mapsto \kappa_0-e_0-1, \quad e_n \mapsto -e_{L-n} \quad (n \neq 0), \\
\kappa_0 \mapsto \kappa_0, \quad \kappa_n \mapsto - \kappa_{L-n}  \quad (n \neq 0), 
\quad \theta \mapsto \kappa_0-\theta,
\\
\displaystyle
q_n \mapsto \frac{s p_{L-n} }{q_{L-n}p_{L-n}-\kappa_{L-n} },
\quad
p_n \mapsto \frac{q_{L-n} (\kappa_{L-n}-q_{L-n}p_{L-n})}{s}.
 \end{array}
\right.
\]
\end{remark}

\allowdisplaybreaks[0]

\appendix
\section{Case $L=2$: the Garnier system}
If $L=2$ then the canonical Hamiltonian system ${\cal H}_{L,N}$
is equivalent to the Garnier system in $N$ variables
(see \cite{gar12}).
This fact is guaranteed by the Lax formalism given in Sect.~\ref{sect:lax}.
But, however, our polynomial Hamiltonian function (\ref{eq:caneq_ham}) is different from that given in \cite{ko84}
(see also \cite{iksy91}).
In this appendix
we describe explicitly the canonical transformation
between the two Hamiltonian systems.

First we concerns the general $(L,N)$ case.
Define the canonical transformation
$(q_n^{(i)},p_n^{(i)},H_i,s_i) \mapsto (Q_n^{(i)},P_n^{(i)},\widetilde{H}_i,s_i)$
by
\begin{align*}
Q_n^{(i)}
&=-s_i \frac{p_n^{(i)}}{p_n^{(0)}}
\left(
=-s_i\frac{b_n^{(i)}}{b_n^{(0)}}
\right),
\\
Q_n^{(i)} P_n^{(i)}
&=-q_n^{(i)}p_n^{(i)}
\left(=b_n^{(i)}c_n^{(i)}
\right),
\\
\widetilde{H}_i
&= H_i -\sum_{n=1}^{L-1} \frac{q_n^{(i)}p_n^{(i)}}{s_i}.
\end{align*}
Clearly the new Hamiltonian function
$\widetilde{H}_i$ becomes again a polynomial in 
$Q_n^{(i)}$ and $P_n^{(i)}$.
This canonical transformation is, in short, 
derived from an interchange
of the roles of $b_n^{(i)}$ and $c_n^{(i)}$ in 
the definition 
(\ref{eq:can})
of the canonical variables.
Note that only if $N=1$ it keeps the form of the
Hamiltonian function unchanged, 
thereby giving rise to a birational symmetry; 
see Remark~\ref{remark:add_sym}.

Next we let $L=2$ and write the variables as 
$(Q_1^{(i)},P_1^{(i)})=(q_i,p_i)$
for $i=1,2,\ldots,N$.
The Hamiltonian function $\widetilde{H}_i$
thus takes the following expression:
\begin{align*} 
s_i(s_i-1)\widetilde{H}_i &
\equiv q_i \left(\kappa_1 +\displaystyle{ \sum_{j} q_j p_j}\right)\left(\kappa_1-\theta_0 + \displaystyle{ \sum_{j} q_j p_j} \right)  +s_ip_i(q_ip_i+\theta_i) 
\\
&\quad
 - \sum_{j(\neq i)} R_{ji}(q_jp_j+\theta_j)q_ip_j
   - \sum_{j(\neq i)} S_{ij}(q_ip_i+\theta_i)q_jp_i 
 \\
&\quad
 - \sum_{j(\neq i)} R_{ij}q_jp_j(q_ip_i+\theta_i) 
   - \sum_{j(\neq i)} R_{ij}q_ip_i(q_jp_j+\theta_j) 
 \\
&\quad
 -(s_i+1)(q_ip_i+\theta_i)q_ip_i 
- (\theta_{N+2} s_i+\theta_{N+1}+1)q_ip_i
\end{align*}
modulo some function in only
${\boldsymbol s}=(s_1,\ldots,s_N)$.
Here we put
$R_{ij} = {s_i(s_j-1)}/{(s_j-s_i)}$, $S_{ij} =  {s_i(s_i-1)}/{(s_i-s_j)}$, 
$\theta_{N+1}=d_{1,0}-d_{1,1}-1/2$,
$\theta_{N+2}=d_{1,1}-d_{0,1}-1/2$,
and
$\kappa_1 = 
\left(\sum_{i=0}^{N+2} \theta_i +1\right)/2$.
The symbols 
$\sum_{j}$ and $\sum_{j(\neq i)}$ stand for the summation 
over $j=1, 2,\ldots,N$ and over $j=1, \ldots, i-1, i+1, \ldots, N$, 
respectively.
This is exactly the (usual) polynomial Hamiltonian function for the Garnier system;
cf. \cite{ko84, tsu06}.

\small
\paragraph{\it Acknowledgement.}
I would like to express my sincere gratitude to
Hidetaka Sakai 
for his helpful comments on 
Hamiltonian structure of the Schlesinger systems
and informing me of the literature \cite{jmms}.
I also greatly appreciate
illuminating discussions with
Yasuhiro Ohta and Takao Suzuki. 
This work was partly conducted
during my stay in
the Issac Newton Institute for Mathematical Sciences 
at the program ``Discrete Integrable Systems" 
(2009).

\begin{center}
Faculty of Mathematics, Kyushu University, 
Fukuoka 819-0395, Japan.
\\
e-mail: tudateru@math.kyushu-u.ac.jp
\end{center}

\end{document}